\documentclass{article}
\usepackage[a4paper,bindingoffset=0.2in,%
            left=1in,right=1in,top=1in,bottom=1in,%
            footskip=.25in]{geometry}
\usepackage[utf8]{inputenc}
\usepackage{url}
\usepackage{amsmath}
\usepackage{graphicx}
\usepackage{float}
\usepackage{array,multirow}
\usepackage{hyperref}
\usepackage{amsfonts}
\usepackage{amsthm}
\usepackage{booktabs}
\usepackage{subcaption}
\usepackage[ruled,vlined]{algorithm2e}
\theoremstyle{definition}

\usepackage{xcolor}

\title{
Passenger-Centric Urban Air Mobility: Fairness Trade-Offs and Operational Efficiency
}

\author{Mehdi Bennaceur  \\
	Uber  \\
	\texttt{mehdib@uber.com}\\
	\and 
	R\'emi Delmas \\
	Uber \\
	\texttt{remid@uber.com}\\
	\and 
	Youssef Hamadi \\
	Uber\\
	\texttt{youssefh@uber.com}
	}
\begin{document}
\newcommand{\todo}[1]{{\color{red}#1}}

\maketitle

\begin{abstract}

Urban Air Mobility (UAM) has the potential to revolutionize transportation. It will exploit the third dimension to help smooth ground traffic in densely populated areas. To be successful, it will require an integrated approach able to balance efficiency and safety while harnessing common resources and information. 
In this work we focus on future urban air-taxi services, and present the first methods and algorithms to efficiently operate air-taxi at scale. Our approach is twofold. 
First, we use a passenger-centric perspective which introduces traveling classes, and information sharing between transport modes to differentiate quality of services.
This helps smooth multimodal journeys and increase passenger satisfaction. Second, we provide a flight routing and recharging solution which minimizes direct operational costs while preserving long term battery life through reduced energy-intense recharging.  
Our methods, which surpass the performance of a general state-of-the-art commercial solver, are also used to gain meaningful insights on the design space of the air-taxi problem, including solutions to hidden fairness issues. 

\end{abstract} 


\section{Introduction}

Urban Air Mobility (UAM) will move people and goods by air within and around dense city regions. This new transportation system will help smooth ground traffic despite increasing population densities. UAM will leverage fleets of electric vertical take-off and landing (eVTOL) aircraft operating from small areas, e.g., rooftop \textit{vertiports}, without significant background noise increase.
Besides well known drone for air-delivery efforts \cite{amazon,MURRAY2020368}, several air-taxi services are under active development \cite{Volocopter2019,ehang2020,whitepapercommunityelevate2020,UberWP2016}. They will pool passengers for up to 30 minutes flights, regularly connecting city centers, important suburbs, and airports. Operating air-taxis will require an organized approach able to balance efficiency and safety. It will also require deep levels of automation to enforce safe air traffic control over relatively short times scales and distances \cite{DBLP:conf/lion/Hamadi20}. 

\begin{figure}[htbp]
    \centering
    \includegraphics[width=1\textwidth]{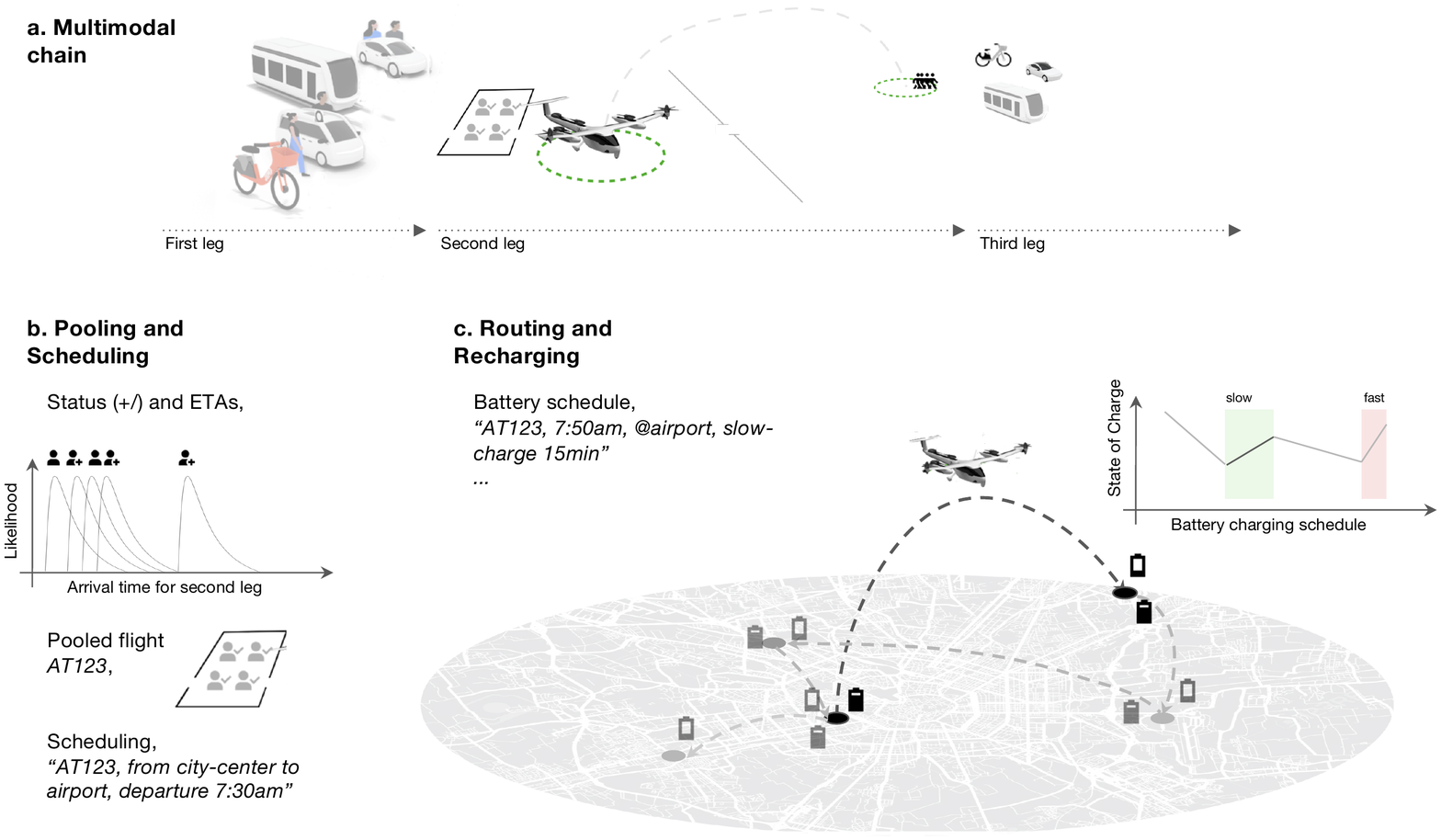}
    \caption{UAM air-taxi. a. Multimodal chain: passengers arrive through some mode, take their flight, then use a third mode to reach their final destination. 
    b. Pooling and Scheduling: flight requests are created grouping compatible passenger demands using expected arrival times (ETAs) and status (class) information.
    c. Routing and Recharging: routing chains different flight requests, ensuring regular battery charging operations.
    }
    \label{fig:abstract}
\end{figure}

Typically, air-taxi flights will be embedded in longer multimodal journeys. Passengers will reach an origin vertiport through some transportation mode, take their urban flight, and complete their journey by leaving the destination vertiport using a last mode. This is presented on top of Figure \ref{fig:abstract}. In this work, our objective is to exploit this rich context to provide a complete solution for air-taxi services. We consider two complementary problems.

The first problem, "Pooling and Scheduling", performs the grouping of passengers demands into common flights, adequately scheduled to meet their expected arrival time at origin vertiports and potential constraints on arrival time at destination. The number of created flights is minimized. This is presented on the bottom left part of Figure \ref{fig:abstract}. We consider this problem through a passenger-centric perspective which puts air-taxi passengers at the centre of future solutions, where their waiting time between their arrival at the vertiport and the effective flight departure is minimized. Passenger-centricity builds on the notion of Quality of Service (QoS) and extends it to include non-strictly operational factors that add to the ability of service providers to differentiate their offers \cite{DBLP:journals/pubtrans/CamachoFRRB16}.
It is valuable to passengers and advantageous for providers, allowing the creation of different service classes associated to different expected waiting times, e.g., Premium vs Regular.

The second problem, "Routing and Recharging", deals with the routing and battery recharging of electrical aircraft. The goal is to make the best use of an eVTOL fleet to serve flights created in the first problem while minimizing underlying operational costs, such as (but not limited to) deadheads. In particular, we pay special attention to the energy rates and duration used for lithium-ion battery charging actions. Batteries could for instance be charged in two modes, \textit{slow} or \textit{fast}. The fast mode is more energy-intense and rapid than the slow one, but is known to have a lower efficiency
and cause accelerated capacity and power fade \cite{tomaszewska2019lithium}. Therefore it is preferable not to use it when possible. Our methods make the best use of a fleet, minimizing operational costs, and preserving long term battery life.

The previous problem decomposition fits envisioned general air-taxi operational workflow. \emph{Flight brokers} will market flights, and collaborate with \emph{fleet operators} able to execute them.
Therefore, brokers need efficient pooling of passenger demands to minimize the number of flights requested to operators, and operators need efficient routing and recharging solutions to efficiently perform their orders. These distinct interests can be associated to different objective functions which have to be optimized separately.

This paper presents the following contributions:
\begin{enumerate}

\item We define, the UAM air-taxi problem. It includes several sub-problems: demand consolidation into joint flights, vehicle selection, routing, and energy management.

\item We provide the first Mixed Integer Linear Programming models for this problem, decomposed into "Pooling and Scheduling" with passenger-centric metrics, and "Routing and Recharging" with battery health considerations.

\item We propose the first dedicated algorithms for solving this problem. A beam-search algorithm for the "Pooling and Scheduling" part, and a variable neighborhood search algorithm for the "Routing and Recharging" part. 

\item We apply and test our methods against foreseen air-taxi infrastructures, surpassing the performances of a state-of-the-art MILP solver, Gurobi 9.1.0. 

\item We provide meaningful insights on the design space of the air-taxi problem, in particular we study the effects of introducing passenger classes. We find how QoS is impacted by class parameters and characterize intra-class fairness issues.

\end{enumerate}

In the following, we present related works in Section 2. Section 3 describes our "Pooling and Scheduling" and "Routing and Recharging" problems in the general Mixed Integer Linear Programming formalism. Algorithms are presented in Section 4 which describes a new beam-search algorithm to efficiently pool passengers in joint flights, and a new variable neighborhood search algorithm to effectively allocate and route eVTOLs given a flight schedule produced by the Pooling and Scheduling phase. Before our conclusion in Section 6, computational results and their analysis are presented in Section 5. They demonstrate the superiority of our algorithms against a state-of-the-art MILP solver, Gurobi 9.1.0, and provide great insights on the air-taxi problem design space.

\section{Related Works}
Although recent research has been conducted on air-taxi \cite{uam-problem-optim,uam-problem-schedule-or-demand,uam-evtol-electricity-market}, it has not been extensively studied and many questions remain unanswered. In particular, the pooling and scheduling part of the problem, the electrical energy used, along with the possibility to use different charging modes (i.e., fast or slow), and the on-demand nature of the service, makes the problem different from traditional vehicle routing problems. 
Therefore, this section presents relevant works from similar domains. We categorize problems according to several aspects : the type of decisions involved, the associated objective functions, and the algorithms. We conclude the section with a discussion on passenger-centric works in transportation.

\subsection{Type of Decisions}

In the air-taxi UAM problem, we are dealing with four main type of decisions, pooling decisions which determine which passengers fly together, scheduling decisions to determine flight dates, routing decisions to allocate eVTOLs to feasible sequences of flights, and battery recharging decisions needed to execute the flights.

Routing decisions are largely present in problems studied in the literature, in particular in vehicle routing problems (VRP), see \cite{vrpsurvey} for a review. In these problems, a fleet of vehicles is used to serve a set of customers under various constraints. There, typical decision variables mainly relate to deciding which customers each vehicle is going to serve and how much time each vehicle is going to spend waiting at each node it visits \cite{vrpsurvey}. Customers are located on nodes in the graph representing the problem and vehicles must visit these nodes to serve the demands, in some cases however (see \cite{thesisann} for examples) demands may be represented on edges rather than nodes. This is the case in arc routing problems (ARP), which are in fact equivalent to VRP \cite{vrptocarp}, where the decision variables remain similar. In these problems, vehicles are initially located at a depot and must end their route there as well. 
A variant is the Pickup and Delivery Problem (PDP) \cite{pdpsurvey} where customers make transportation requests, each with a pickup and a delivery node. In the context of a taxi service, it is more commonly called Dial-a-Ride Problem (DARP). VRPs can be viewed as a special case of PDP in which pickup nodes are customers nodes and all delivery nodes are located at the depot.

In the late 2000s, the DayJet corporation operated an on-demand regional jet taxi service in the US. The company allowed customers to book seats for regional flights operated through a fleet of six-seater jets. With this new service, new optimization problems appeared and shed light on the Dial-a-Flight (DAF) problem \cite{dayjet1,dayjet2,daf}. In \cite{dayjet1}, the problem is modeled by a graph in which nodes represent elementary actions (e.g., loading passengers, standing-by, etc.), the main decision variables assign arcs to aircraft, constructing a path for each of them. Another approach \cite{daf}, makes use of \textit{composite variables} which represent sub-routes of aircraft and contain more information than elementary variables.

When handling a fleet of vehicles, another important factor is the energy management of the vehicles, be it fuel or electricity. From this arises a second type of decisions which deal with battery charge or fuel level management. This is traditionally not taken into account in most VRP where fuel load is typically enough to carry out a full day of operations. \cite{daf}, for example, studies a DAF problem without integrating refueling as a decision and assumes that aircraft can always carry enough fuel for all flights of the day. Others however include these decisions, \cite{bousonfuel} proposes a heuristic to solve a VRP with time windows (VRPTW) while making refueling decisions on both quantity and timing.

During the last decade, research has been conducted on electric versions of vehicle routing problems which might require several recharges during operations \cite{electricvrptw,electricvrp2}. For example, \cite{electricvrptw} introduces a VRP with electric vehicles and recharging stations, thereby considering new decisions and making the optimization more complex. The authors embed the problem within a graph containing both customer nodes and recharging nodes, then use a primary decision variable which assigns arcs to vehicles, making the recharging decision implicit and activated when a vehicle visits a recharging node. \cite{bruglieri2015variable} also includes the charging times as decision variable in an electric VRP. In a UAM context, recharging decisions are also present in \cite{uam-evtol-electricity-market,uam-problem-optim} and \cite{uam-problem-schedule-or-demand} in which eVTOLs can visit recharging pads or traverse specific arcs modeling recharging. \cite{uam-evtol-electricity-market} introduces a new energy-related decision which is to allow idle eVTOLs to transfer power to the electricity grid, acting as a regulation tool. Authors also define a price decision variable which in turn determines the demand quantity through a demand model.

An important characteristic of the air-taxi UAM problem is its on-demand and multimodal nature, which implies a critical pooling phase. It deals with grouping passengers together into flights and resembles Bin Packing problems with conflicts (BPC) \cite{bin-packing-conflict}, in which items have to be assigned to a minimum number of bins of limited capacity, some items excluding  others in a bin. The passenger pooling problem decides which passengers to group together, under aircraft capacity constraints and incompatibility constraints based on waiting time.

In Pickup and Delivery or Dial-a-Ride problems \cite{pdpsurvey,dar}, pooling naturally emerges as a by-product of the optimization, but the generated groups are not subject to specific constraints or optimization objectives.
In an air-taxi UAM context however, pooling passengers is a primary objective for economic reasons and needs to be finely controlled to provide the best quality of service to passengers. Additionally, contrary to the DAR setting, pickup and delivery locations are not highly dynamic and unknown. Under this light, the air-taxi UAM problem is closer to the DAF problem. \cite{daf} studies the DAF problem restricted to a specific case study in which most passenger bookings could be aggregated in one only way and a brute force search could be used for the few remaining cases. While \cite{dayjet1} considers a more complex online pooling problem, in which passengers are informed in real time if their demand can be served and at what price, they do not disclose their methodology. \cite{uam-problem-optim} studies an air-taxi problem and allows for pooling during optimization but is not constrained by the fleet size and considers only one vertiport, not capturing interactions between vertiports in a network.

\subsection{Objective Functions}

In routing problems, objective functions are usually centered around minimizing the actual cost of implementing the solution, subject to a set of operational constraints \cite{vrpsurvey}. 
When dealing with passenger transportation, authors usually tend to favor constraints and objectives that are more service oriented. These may be related to waiting time or hard time windows preventing deviation from the required time of pickup. \cite{darpreview} gives an overview of DARP and describes the modeling of constraints and objectives related to managing user inconvenience. For example, \cite{patienttransport} studies an application of DARP and, to account for service quality, authors define an objective function which is a weighted sum of three metrics related to service quality. 

The air-taxi pooling and scheduling problem we present occurs at the meeting points of two legs, i.e., two different modes of transportation - see Figure \ref{fig:abstract}. In this context, the issue of potential passenger delays emerges. This point has been the focus of some recent research, for example \cite{kucharski2020if} studies the effect of passenger lateness on individual and overall delays within a shared ride. The authors consider a fixed scenario where a sequence of pickup and deliveries is known and evaluate the impact of passenger delays, depending on their position within the trip.  \cite{heilporn2011integer} considers a variant of DARP where customers might not show up on time with stochastic delays, this induces a delay cost which is included in the objective function to be minimized. In a UAM context, \cite{uam-problem-optim} adopts a passenger-centric approach by minimizing the passenger waiting time at the vertiport, however the authors do not consider fleet operating costs.

The routing part of the problem involves different costs which can be operational or related to the batteries. \cite{uam-evtol-electricity-market} uses a cost-based objective function, computing the cost of running an air-taxi service, including operating and electricity-related costs and revenues. Authors consider using the eVTOL fleet as a mean to help frequency regulation of the electric grid, thereby generating additional revenues. An additional source of complexity comes from battery management. 
In particular, \cite{zhang2006effect} studies various charging protocols for Li-ion batteries, and shows that some can have a significant impact on their health. These charging protocols can be optimized to achieve shorter charging times \cite{dung2010ilp} but also to achieve a trade off between charging efficiency and battery state of health (SoH) \cite{liu2018charging}. For example, \cite{barco2013optimal} includes a battery degradation cost in a vehicle routing problem applied to an airport shuttle service using electric vehicles, but without including different charging modes. The authors first compute a set of optimal routes with respect to energy consumption, then determine route to vehicle assignments and vehicle's charging schedules. The charging schedules are chosen to minimize an objective including a battery degradation cost which is defined as a function of their temperature, states of charge and discharge. 

In our problem, charging decisions have to be made jointly with routing decisions, including the possibility of different charging modes (fast vs slow), while considering both routing and charging costs and impacts on battery health. Furthermore, the use of aircraft comes with several aviation-specific constraints, for example only part of the battery capacity will be considered as actually usable, forcing more frequent recharge operations. Furthermore, eVTOL batteries are still under development and will not be the same as those used in traditional vehicles, precluding any direct reuse of existing battery models, such as \cite{barco2013optimal}.

\subsection{Algorithms}

The problems presented in this paper are combinatorial by nature. Such problems are often formulated using mixed integers linear programs (MILP). These formulations allow the use of state-of-the art general purpose solvers, e.g. Gurobi, CPLEX. While providing optimal solutions, these methods typically require important computational resources, \cite{thesisann} exposes these limits while studying CARPs. To overcome these limitations, Column Generation algorithms have been successfully applied to routing problems \cite{CGhvrp}.

While not providing guarantees on optimality, heuristics and metaheuristics are incomplete methods which often achieve good quality versus time trade off and have been found successful to solve routing problems \cite{thesisann}. 
Heuristic algorithms are designed solely to solve a specific problem, providing good solutions quickly. \cite{bousonfuel} presents a heuristic algorithm to solve the VRPTW which is then coupled with an exact method to integrate refueling decisions. 
Metaheuristics are more general and less problem-specific and can therefore be used to solve different optimization problems. These methods are also anytime, yielding valid solutions while running, this feature can be desirable when time budgets are very low. In particular, Variable Neighborhood Search (VNS) \cite{affi2018variable,hansen2017variable} is a state-of-the-art algorithm for many optimization problems and has been successfully applied to VRP. In \cite{electricvrptw}, a Variable Neighborhood Search is coupled with Tabu Search (TS) to solve an electric VRPTW, the authors report good performances as well as a positive impact of mixing both VNS and TS.

More recently, learning algorithms have been used to address routing problems.  \cite{mdpvrp} uses reinforcement learning to address a variant of the VRP in which customer requests may appear randomly during service. By framing the problem as a Markov decision process (MDP), the authors naturally embed the uncertainty in the model. Others, such as \cite{learning-carp} propose to use Graph Neural Networks to first compute an embedding for the Graph representing the problem and then use it as input to a policy gradient algorithm used to solve a CARP. While these methods open interesting perspectives (e.g., generalization, faster solving times), they typically require a large number of examples, very large computational resources during training and more research is needed to achieved the desired generalization.

\subsection{Passenger-Centric Approaches}

Passenger preferences have been taken into account recently in some transport applications. \cite{passenger-centric-high-levelview} tries to set the scene and to describe the challenges of placing passengers at the center of transportation systems, identifying four main axes: passenger experience, safety and security, multimodality, and mobility as a service. Taking into account passenger preferences in transportation systems requires understanding several aspect of the service. For example, recent research \cite{alonso2020value} has been conducted to evaluate the passengers' value of time and reliability in a pooled on demand urban ride sharing service.

\cite{passenger-centric-train} studies the trade-off between operator profit and passenger utility in a railroad network, finding that significant improvement in passengers satisfaction can be achieved without dramatic impact on profits. The authors study cyclic timetable (i.e., fixed, repeating schedules) and non-cyclic (i.e., non-repeating schedules) mentioning that cyclic schedules have the advantage to be easy to remember and may provide an utility to passengers. In the context of UAM, it has been studied in \cite{uam-problem-schedule-or-demand} whether the service of urban air taxis should be on demand or scheduled, but not accounting for heterogeneous passenger preferences. In traditional aviation, some works adopt a passenger-centric approach to optimize air traffic flow management decisions such as Ground Holding Policies (GHP). \cite{flight-flow-management} evaluates the relevance and effect of using different objective functions on passenger or flight centered metrics. The authors discuss the consequences of using passenger-centered objectives and conclude that it might lead to unfair management between different airlines. 
\cite{passenger-centric-ground-holding-atfm} addresses fairness concern by proposing an objective function which minimizes passengers waiting times and delay costs while ensuring the computed GHP is not deviating too much from the aviation traditional first scheduled first served (FSFS) basis.
\newline

Our problem requires an important pooling and scheduling phase to consolidate individual demands which has to be finely controlled to ensure a good quality of service, differing from PDP and DARP \cite{pdpsurvey,dar}. In the second phase, routing and recharging decisions have to be made jointly and considering both operational costs and battery health. Additionally, battery usage is restricted due to aviation-specific safety constraints and models previously developed for other problems \cite{barco2013optimal} cannot directly be reused. Recent works related to air taxis \cite{uam-problem-schedule-or-demand,uam-evtol-electricity-market,uam-problem-optim} do not fully capture these aspects. In the following sections, we formalize our problem and present dedicated algorithms to tackle it.

\section{Problem Description}

As presented in the introduction, our problem is decomposed into two parts: Pooling and Scheduling followed by Routing and Recharging. This fits envisioned services where flights will be marketed by air-taxi brokers, while operations will be executed by air-taxi operators. Each of them having different incentives for optimization, and only willing to cooperate to some extent. 

In the first part, we deal with the synchronization of different transportation modes. Air-taxi services will often be part of a longer multimodal chain. This is presented in the top part of Figure \ref{fig:abstract}. A first leg is performed through some transportation mode in order to reach an origin vertiport. The second leg is the actual air-taxi flight to some destination vertiport. Eventually, a third leg connects the arrival vertiport to some final destination.
Synchronizing these legs involves grouping compatible passengers together and scheduling appropriate flight departure times. 
It requires information sharing between modes to inform about expected arrival times (ETAs) at origin vertiport. From the user standpoint, a trip from A to B is requested to the broker which then evaluates the feasibility of an air-leg, before proposing the option for booking. The underlying optimization problems can be solved by batching demands or incrementally as demands come in.

The first leg induces the possibility of delays. We address this by scheduling flights such that all passengers arrive on time for boarding with high probability. The first objective of the Pooling and Scheduling part is the maximisation of pooling. To provide users with a good experience we minimize, as a secondary objective, their waiting times based on their status or ticket class, \emph{Regular} or \emph{Premium}. We call \textbf{demands} the transportation demands made by the users to the broker, these can be made until the day before the trip. This is depicted in the bottom left part of Figure \ref{fig:abstract}. We formally define the pooling and scheduling problem in Section \ref{sec:pooling}.

In the second part, a fleet of electric aircraft (eVTOLs) is used to serve the flight requests created in the first problem. We call \textbf{requests} the flight requests made by the broker to the operator, these requests are the result of the consolidation of demands made by the broker in the Pooling and Scheduling phase.
The use of electric aircraft adds a layer of complexity to the standard routing problem. 
In particular, the charging schedule of batteries can largely impact the quality of the routing, the recharging efficiency, and the long term battery health \cite{zhang2006effect}. 
In this paper, we will assume that two charging modes are available at the vertiports : \emph{slow} and \emph{fast}. Using the fast mode allows to fill the battery rapidly but is also associated with higher battery health degradation \cite{tomaszewska2019lithium}. Therefore our objective will include minimizing the number of fast charges used in the solutions. 
This is displayed in the bottom right part of Figure \ref{fig:abstract} and formalized in Section \ref{sec:routingandcharging}.

After this general introduction, we formalize and present our optimization problems in the following. 

\subsection{Pooling and Scheduling}\label{sec:pooling}
This section deals with the pooling of individual demands and their associated flight schedule. As mentioned in the general description users exact arrival time at their origin vertiport is unknown, and we will schedule flights such that all passengers arrive on time for boarding with high probability. Subject to capacity and waiting time constraints, we first minimize the number of created flights (i.e., maximize pooling) and then minimize the overall expected passengers waiting time.

\subsubsection{Assumptions and Notations}

A demand $i$, made by a single user, is characterized by a set of attributes, described in Table \ref{tab:notations-demand}. In addition to these attributes, each demand is related to a pair of origin and destination vertiports. In this section, we consider a set $\mathcal{D}$ of $D$ demands which are all related to the same origin and destination vertiports.

\begin{table}[htbp]
\begin{center}
\begin{tabular}{ c p{9cm} }
 \hline 
 Attribute & Description \\
 \hline
  $w_i$ &  Number of passengers \\ 
  $\nu_i$ &  Arrival time at origin represented by a distribution \\ 
  $s_i$ &  Latest arrival time at destination, -1 if not specified  \\ 
  $c_i$ &  Class (e.g.,  Regular or Premium)  \\
 \hline
\end{tabular}
\caption{\label{tab:notations-demand} Attributes of demand $i$. }
\end{center}
\end{table}

Given these demands, we aim at creating and scheduling a minimum number of flights, such that every demand is assigned to a flight while taking into account passengers waiting times.

\paragraph{Transportation modes synchronization}

We want to create and schedule flights for times at which all passengers assigned to it are present at the origin vertiport with high probability $1 - \delta$. Therefore, any demand $i \in \mathcal{D}$ can only be assigned to a flight which departure time corresponds to the quantile of order $1 - \delta$, noted $q_{1-\delta}^{\nu_i}$, of the distribution $\nu_i$ or later. For any created group of demands $g \subseteq \mathcal{D}$, the corresponding flight's scheduled departure time $f_g$ will thus be such that:

\begin{equation}
    f_g \geq \max_{i \in g} q_{1-\delta}^{\nu_i}
\end{equation}

Demands may have an associated upper bound on arrival time at destination vertiport noted $s_i$ (e.g., users going to an airport), which will forbid assigning them to a flight reaching destination later than $s_i$. The flight time between vertiports being known, this constraint can be translated into a latest departure time noted $t_i^{\mathit{max}}$.

\paragraph{Quality of service}
In traditional airlines, passengers can pay a premium to have a better flight experience including more space or better catering. 
In the air-taxi setting however, restricted eVTOLs size and design will limit this kind of separation between passengers. 
Therefore, we relate the premium to the \textit{waiting time} rather than comfort. A premium demand will be given stronger guarantees on its waiting time at the origin vertiport. More precisely we will require that premium demands have an expected waiting time of at most $t_{\mathit{premium}}$ minutes, while for regular demands it may go up to $t_{\mathit{regular}}$ minutes with $t_{\mathit{premium}} < t_{\mathit{regular}}$. The expected waiting time for demand $i \in g$ can be expressed as:
\begin{equation}
    WT_i^g = f_g - \mathbb{E}[\nu_i]
\end{equation}

Furthermore, to account for heterogeneous values of time (see \cite{alonso2020value}), the waiting time of demand $i$ is weighted by a factor $\alpha_{c_i}$ (see Equation (3)) depending on its premium or regular class. Adjusting these weights determines how the overall waiting time is distributed between premium and regular demands (see Section 5.2).

\begin{figure}[htbp]
    \centering
    \includegraphics[width=0.8\textwidth]{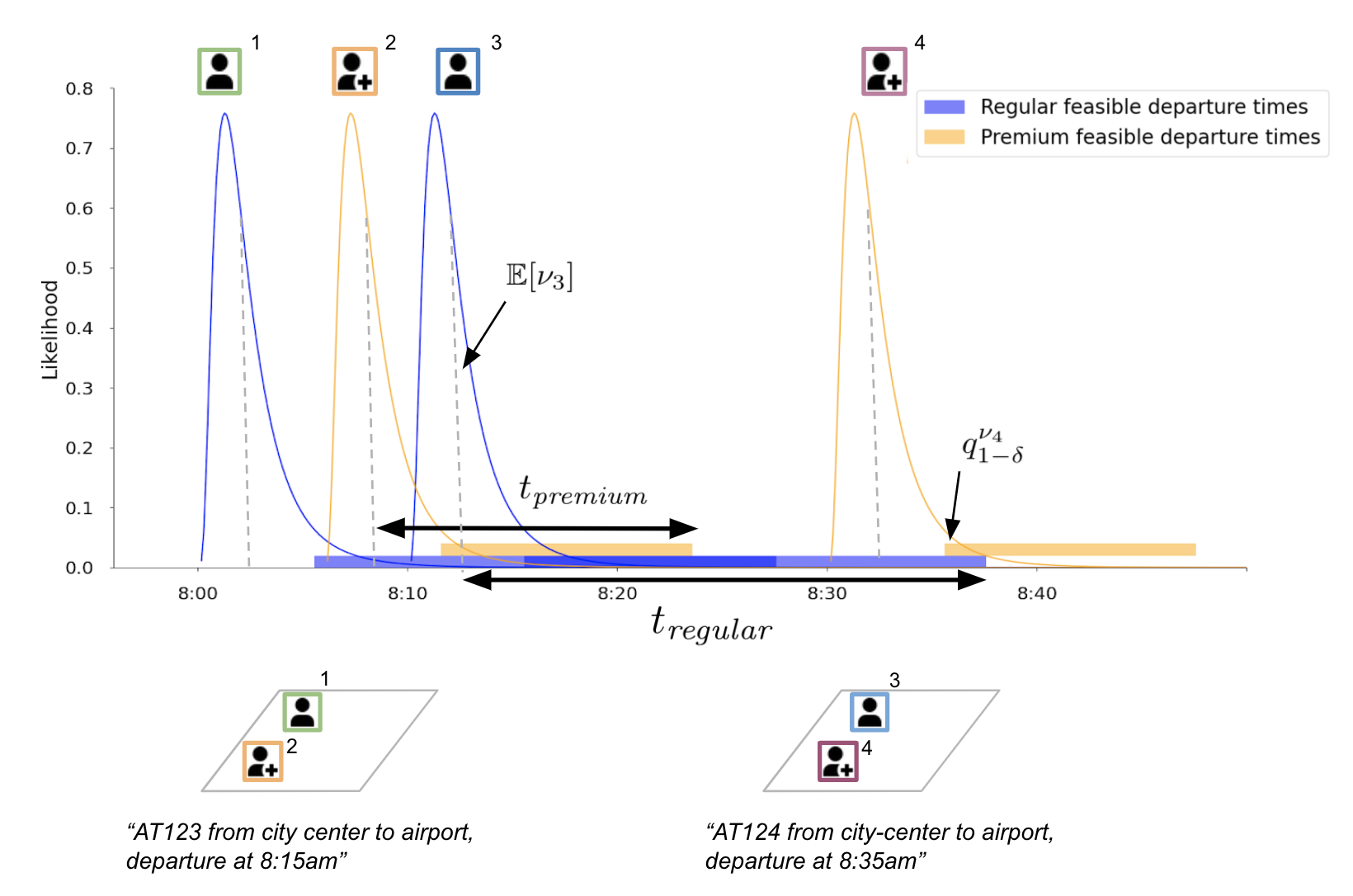}
    \caption{Example for Pooling and Scheduling. A set of 4 demands, 2 premiums (+) and 2 regulars, is pooled into 2 flights.}
    \label{fig:pooling-constraints}
\end{figure}

\hfill

An example of pooling and scheduling is given in Figure \ref{fig:pooling-constraints}. Four demands are represented with their respective arrival time distribution, 2 premiums (+) and 2 regulars. Maximum expected waiting times $t_{\mathit{regular}}$ and $t_{\mathit{premium}}$ are represented for demand 2 and 3, highlighting a tighter constraint for the premium one. Colored zones, deduced from (1) and (2), close to the x-axis represent possible departure time for each demand, yielding the possible grouping combinations. Notice that demand 3 could be pooled with demand 1 and 2, however this would result in more waiting time for 2, which is premium. It is therefore preferable to pool 3 with 4.  

\subsubsection{Model}\label{sec:model-pooling}

We model the pooling problem as a variant of the Bin Packing Problem with Conflicts \cite{bin-packing-conflict}. In our case, bins are the groups created (i.e., containing passengers that will fly together) and bin capacity corresponds to aircraft's capacity.
Demands classes and latest departure times create a set of exclusion constraints preventing some demands from being pooled together.

Bins are indexed by $k$, with as many bins as individual demands, the binary variable $y_k$ indicates the usage of bin $k$ and $f_k$ is its associated flight departure time. C denotes the capacity of one aircraft. Our decision variables will target assignment choices and flight departure times. Binary variable $x_i^k$ will be equal to 1 if demand $i$ is assigned to bin $k$, and the real-valued variable $f_k$ is the flight departure time for bin $k$. We define the overall weighted expected waiting time ($OWT$) to be the sum of expected waiting time over all demands weighted by a factor $\alpha_{c_i}$ depending on demand's class.
\begin{equation*}
     OWT = \sum_{k=1}^{D} \sum_{i=1}^{D} x_i^k \cdot (f_k - \mathbb{E}[\nu_i]) \cdot \alpha_{c_i}
\end{equation*}

Our first objective is to minimize the number of bins used and the second to minimize the OWT. We encode this lexicographic optimization into a single objective function (3) using stratification, $\lambda_p$ being the penalty parameter for using one bin.

\begin{equation}
   \min \sum_{k=1}^{D} y_k \cdot \lambda_p + \sum_{k=1}^{D} \sum_{i=1}^{D} x_i^k \cdot (f_k - \mathbb{E}[\nu_i]) \cdot \alpha_{c_i}
\end{equation}

\textbf{Constraints :}
\begin{align}
     \sum_{i = 1}^{D} x_{i}^k \cdot w_i \leq C \cdot y_k \;, \; \forall k \in \{1..D\} \\
     \sum_{k=1}^{D} x_{i}^k = 1 \;, \; \forall i \in \{1..D\}\\
    x_i^k \cdot f_k \leq t_i^{max} \; , \; \forall i, k \in \{1..D\}\\
    x_i^k \cdot (f_k - \mathbb{E}[\nu_i]) \leq t_{c_i} \; , \; \forall i, k \in \{1..D\}\\
    f_k \geq q_{1-\delta}^{\nu_i} \cdot x_i^k \; , \; \forall i, k\in \{1..D\} \\
     y_{k+1} \leq y_k \;, \; \forall k \in \{1..D - 1\}\\
     f_{k+1} \leq f_k \;, \; \forall k \in \{1..D - 1\}
\end{align}
\begin{align*}
     y_k, x_{i}^k  \in \{0, 1\}, f_k \in \mathbb{R}
\end{align*}
    
The pooling and scheduling objective function (3) contains two components, the first 
one penalizes group creation and the second corresponds to the overall weighted expected waiting time. The capacity constraint is enforced through (4) which also monitors group usage and (5) ensures each demand is assigned to exactly one group. Demands which have a latest possible arrival time may not be assigned to a flight departing too late (6). A maximum expected waiting time is enforced by constraint (7) where $t_{c_i}$ represents the threshold for class $c_i$. Conditions for multimodal synchronizing derived in (1) are enforced in (8) so that all passengers are likely enough to be all present for boarding at the scheduled date. Finally, as group indexing is meaningless many symmetries are present in the problem (i.e., any index permutation), we break them using constraints (9) and (10). 

For efficiency, expectations and quantiles present in the formulation are computed in advance and stored in a lookup table. The objective function (3) contains a non-linear term due to the product between $x_i^k$ and $f_k$, however since the first is a binary variable and the second a continuous variable the product can be linearized by introducing an additional variable and using the big-M method. The same method is used to linearize constraints (6) and (7).

\hfill

Once the pooling and scheduling problem is solved, we have a set of groups with their associated departure times. To each of these groups, we associate a \textit{request} which is characterized by origin and destination vertiports and flight departure time.
These requests are the input to the Routing and Recharging problem presented in the following sub-section.

\subsection{Routing and Recharging}
\label{sec:routingandcharging}

This section deals with the routing and recharging of aircraft for a given set of requests obtained after the pooling and scheduling step presented in Section \ref{sec:pooling}. 
Our goal is to maximize the number of flights requests served first, scheduling charging operations while minimizing the number of fast charge used second and minimize the operational cost third. Similar to the previous section, we embed this lexicographic optimization into a single objective function.

\subsubsection{Assumptions and Notations}
\paragraph{Batteries} \label{sec:batteries}

\begin{figure}[htbp]
    \centering
    \includegraphics[width=\textwidth]{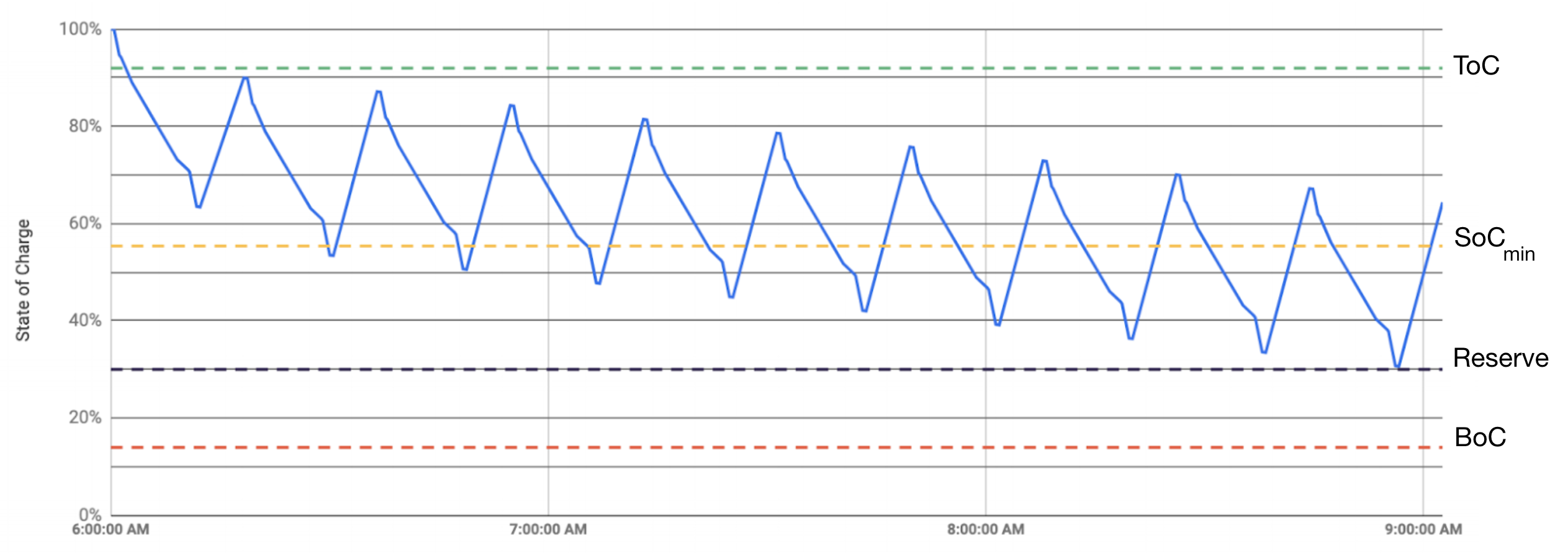}
    \caption{State of Charges constraints \cite{uam-problem-optim}}
    \label{fig:soc}
\end{figure}

The recharging of lithium-ion batteries will be constrained by both manufacturer instructions and safety considerations. We illustrate the different State of Charge (SoC) constraints in Figure \ref{fig:soc}.
Batteries charging should stop at a certain ceiling above which charging rate decreases, we will refer to this value as the Top of Charge (ToC). The battery's Bottom of Charge (BoC) corresponds to a level under which behavior is unreliable and should not be considered as usable energy. These two values can be used to define the acceptable domain of variables measuring aircraft's State of Charge (SoC). To ensure safe flights, a minimum value $SoC_{min}$ is required for takeoff and a minimum reserve should be kept at all times for contingency management. We assume that the upper-bound energy consumption of traversing each route is known in advance and we therefore can define a conservative $SoC_{min}$ which ensures the minimum reserve is never used.

Furthermore, we assume that smart chargers may be used in two modes: \textit{fast} or \textit{slow}. The fast mode is more intense and rapid than the slow one, but also degrades the battery faster \cite{tomaszewska2019lithium}, therefore it is preferable not to use it when possible. The recharging is assumed to increase the SoC linearly with time, with a slope depending on the charge mode.

The state of charge after a flight between two vertiports $i$ and $j$ is a complex non-linear function of the state of charge at takeoff, flight duration and speed, and other factors such as battery temperature, weather conditions and aircraft load. In this paper, we will use a conservative upper bound noted $g_{i \xrightarrow{} j}$ for this quantity.

While our battery charge model is simple, it can be easily refined with more detailed  aviation-specific battery models (piecewise-constant or piecewise-linear) when they become available.

\paragraph{Routes}\label{sec:routes}

The UAM infrastructure will exploit a set of vertiports $N$. Only a subset $E$ of all the possible routes $N \times N$ will be available for booking. We assume that each landing costs a landing fee $l_i$, dependent on the vertiport. Aircraft are not allowed to deadhead (i.e., to fly empty) unless for re-positioning in order to serve a demand originating at the destination of the deadhead. This assumption allows to focus on decision variables that target specific actions (i.e., serving a demand, recharging) as presented in the following section. The fleet of eVTOLs used will be homogeneous with respect to their characteristics, but all aircraft of the fleet do not necessarily start the day at the same vertiport. We note $\mathcal{R}$ the set of requests created after the pooling and scheduling process (Section \ref{sec:pooling}). Each $r$ in $\mathcal{R}$ is characterized by a location of departure, arrival, a departure time and arrival time, i.e., $r = (r^-, r^+, t_{r^-},  t_{r^+})$, with $(r^-, r^+) \in E$.
The operating cost of an aircraft is assumed to be proportional to the amount of time it is being used (i.e., flown) and each request $r$ will have a monetary value of $\beta_r$.
To ensure there is always enough time for passengers to ingress, egress and for the pilot to perform necessary pre-flight checks, we define $\Delta$, a minimum time spent on the ground between two flights for all aircraft. For completeness, notations and their description are presented in Table \ref{tab:notations}.

\begin{table}[htbp]
\begin{center}
\begin{tabular}{ c p{12cm} }
 \hline 
 Notation & Description \\
 \hline
  $T$ & Sequence of one-minute time steps over one service day \\ 
  $N$ &  Set of vertiports \\ 
  $E$ &  Set of legs.\\ 
  $V$ &  Set of aircraft (eVTOLs) \\ 
  $\mathcal{R}$ &  Set of requested flights \\ 
  $\eta$ &  Cost in \$ per minute of operation for one eVTOL  \\ 
  $l_i$ &  Cost in \$ for one aircraft to land at $i \in N$ \\ 
  $\beta_r$ & Gain in \$ of serving request $r$  \\ 
  $t_{i \xrightarrow{} j}$ & Time in minutes to fly from vertiport $i$ to $j$  \\ 
  $g_{i \xrightarrow{} j}$ & SoC units consumed to fly from $i$ to $j$ (upper bound)\\
  $\Delta$ & Minimum time required between a landing and a takeoff, in minutes  \\
  $SoC_{min}$ & Minimum SoC needed for an eVTOL to take off  \\
  $ToC$ & Maximum SoC allowed to preserve batteries  \\
  $p_e$ & Price in \$ of charging 1 SoC unit.\\
  $\gamma_s$ & Charging rate in slow mode (in SoC/minute)\\
  $\gamma_f$ & Charging rate in fast mode (in SoC/minute)\\
 \hline
\end{tabular}
\caption{\label{tab:notations} Notations }
\end{center}
\end{table}

\subsubsection{Model}

\begin{figure}[htbp]
    \centering
    \includegraphics[width=1\textwidth]{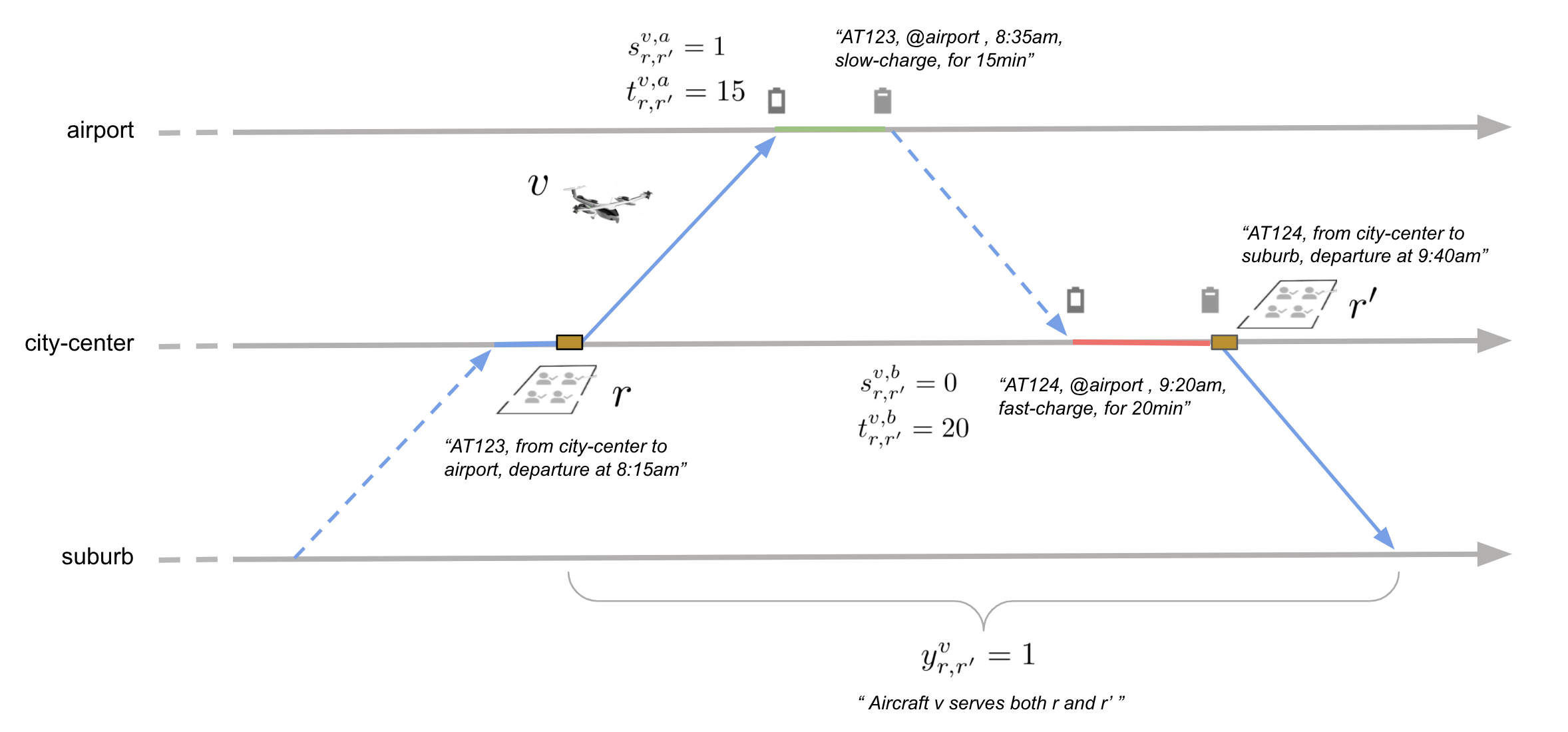}
    \caption{Example for Routing and Recharging. Two flight requests $r$ and $r'$ are served by the same eVTOL $v$.}
    \label{fig:routing-variables}
\end{figure}

We propose a MILP formulation in which decisions are variables associated with pairs of requests $(r, r')$ and target service actions (i.e., serving a request), charging mode and duration selection (see in Table \ref{tab:dvar}). For clarity, decision variables and their role are illustrated in Figure \ref{fig:routing-variables} where an aircraft $v$ is represented serving two requests $r$ and $r'$ with two recharges, dashed arrows representing deadheads.

\begin{table}[htbp]
\begin{center}
\begin{tabular}{ c p{12cm} }
 \hline 
 Variable & Description \\
 \hline
  $y_{r, r'}^v$ & Binary variable indicating that aircraft $v \in V$ sequentially serves request $r$ and $r'$ in that order. \\ 
  $s_{r, r'}^{v,a}$ & Binary variable indicating that aircraft $v \in V$ will use a slow charging mode after serving request $r \in \mathcal{R}$, when connecting $r$ to $r'$. If a slow charge is not used, then it is fast. \\ 
  $s_{r, r'}^{v,b}$ & Binary variable indicating that aircraft $v \in V$ will use a slow charging mode before serving request $r' \in \mathcal{R}$, when connecting $r$ to $r'$. If a slow charge is not used, then it is fast. \\
  $t_{r,r'}^{v,a}$ &  Continuous variable indicating the time, in minutes, aircraft $v$ charges after serving request $r \in \mathcal{R}$, when connecting $r$ to $r'$. \\ 
  $t_{r,r'}^{v,b}$ &  Continuous variable indicating the time, in minutes, aircraft $v$ charges before serving request $r' \in \mathcal{R}$, when connecting $r$ to $r'$. \\
 \hline
\end{tabular}
\caption{\label{tab:dvar} Decision Variables}
\end{center}
\end{table}

For the variables to be coherent, we can augment the set $\mathcal{R}$ with fictitious requests, noted $U$, associated with the departure location of each aircraft $v \in V$ : $r_{0}^v = (start_v, start_v, t_{start}^v, t_{start}^v)$. Departure locations are given as problem input for every aircraft. Aircraft's batteries may be charged after each landing using either the slow or the fast mode. Since arbitrary deadheads are not considered (see Section \ref{sec:routes}),
charging after landing means charging \textit{after} or \textit{before} serving a request. Therefore, we use two continuous variables, $t_{r,r'}^{v,a}$ and $t_{r,r'}^{v,b}$ to model the charging times of aircraft: $t_{r,r'}^{v,a}$ is the charging time after serving r, spent at its destination vertiport, and $t_{r,r'}^{v,b}$ is the charging time before serving r', spent at its source vertiport. Similarly, two binary variables $s_{r, r'}^{v,a}$ and $s_{r, r'}^{v,b}$ indicate which charging mode is utilized after $r$ and before $r'$. These variables jointly determine the quantity (in SoC units) of electricity purchased for recharging aircraft $v$ between serving $r$ and $r'$, noted $b_{r, r'}^v$. Additionally, when deadheads are required we assume that they happen as soon as possible (see Figure \ref{fig:routing-variables}), that is to say after at least $\Delta$ minutes and after the chosen charging time is completed.

The cost of operating aircraft being known in advance, we can compute the cost $c_{r, r'}$ for an aircraft to serve $r'$ when the last served request is $r$, i.e., cost of \textit{activating} $y_{r, r'}^v$.

\begin{flalign*}
    c_{r, r'} &= \begin{cases}
               \infty, \quad if \quad t_{r^+} + t_{r^+ \xrightarrow{} r'^-} + 2 \Delta \cdot 1_{r^+ \neq r'^-} + \Delta \cdot 1_{r^+ = r'^-} > t_{r'^-} \\
               \eta \cdot t_{r'^- \xrightarrow{} r'^+} + l_{r'^+} - \beta_{r'}, \quad if \quad r^+ = r'^- \bigwedge t_{r^+} + \Delta \leq t_{r'^-} \\
               \eta \cdot (t_{r^+ \xrightarrow{} r'^-} + t_{r'^- \xrightarrow{} r'^+}) + l_{r^+} + l_{r'^+} - \beta_{r'}, \quad if \quad r^+ \neq r'^-  \bigwedge t_{r^+} + t_{r^+ \xrightarrow{} r'^-} + 2 \Delta \leq t_{r'^-}
            \end{cases}
\end{flalign*}

By using an infinite cost when there is not enough time to connect $r$ to $r'$ while respecting the minimum amount of time between successive landing and takeoff $\Delta$, we ensure no \textit{invalid} connection is made. For valid connections, the cost is the sum of landing fees and operating costs of the route starting from $r^+$ (i.e., after the service of $r$) and ending at $r'^+$ (i.e., after the service of $r'$), minus the value $\beta_{r'}$ associated with $r'$.

As previously mentioned the first objective is to minimize the number $n_{u}$ of unserved requests. Furthermore, the use of a fast charging rate having a stronger degradation effect on batteries (cf. \cite{tomaszewska2019lithium}), we want to minimize the number $n_f$ of fast charges, in order to maximize the lifespan of eVTOL batteries. We associate penalty coefficients $\alpha_{f}$ and $\alpha_{u}$ to $n_f$ and $n_u$ in the objective function (11). We perform the following lexicographic optimization by setting the penalty weights appropriately: the number of served demands is maximized, the number of fast charges is minimized, the operating and charging cost are minimized.

\begin{equation}
    \min \quad n_{u} \alpha_{u} + n_f \alpha_{f} + \sum_{v \in V} \sum_{(r, r') \in \mathcal{R} \times \mathcal{R}} y_{r, r'}^v \cdot c_{r, r'} + b_{r, r'}^{v} \cdot p_e  
\end{equation}

\textbf{Constraints}

\begin{align}
   \sum_{v \in V} y_{r', r}^v \leq 1, \quad \forall r, r' \in \mathcal{R}\\
   \sum_{r \in \mathcal{R}} y_{r, r'}^v \leq 1, \quad  \forall r', r \in \mathcal{R}, v \in V \\
   \sum_{r' \in \mathcal{R}} y_{r, r'}^v \leq 1, \quad  \forall r', r \in \mathcal{R}, v \in V \\
   y_{r, r'}^v \leq \sum_{p \in \mathcal{R}} y_{p, r}^v, \quad \forall r, r' \in \mathcal{R}, v \in V, r \notin U 
\end{align}
   
\begin{align}
    t_{r, r'}^{v, a} \leq \max(t_{r'^-} - t_{r^+} - t_{r^+ \xrightarrow{} r'^-} - \Delta \cdot 1_{r^+ \neq r'^-}, 0), \; \forall v \in V, \forall r, r' \in \mathcal{R} \\
    t_{r, r'}^{v, b} \leq \max(t_{r'^-} - t_{r^+} - t_{r^+ \xrightarrow{} r'^-} - \Delta \cdot 1_{r^+ \neq r'^-}, 0), \; \forall v \in V, \forall r, r' \in \mathcal{R} \\
    t_{r, r'}^{v, a} + t_{r, r'}^{v, b} \leq \max(t_{r'^-} - t_{r^+} - t_{r^+ \xrightarrow{} r'^-}, 0), \; \forall v \in V, \forall r, r' \in \mathcal{R} 
\end{align}

\begin{align}
    b_{r,r'}^{a, v} = \min(ToC - e_{r^+}^v, t_{r, r'}^{v, a} (\gamma_s s_{r,r'}^{a, v} + \gamma_f (1-s_{r,r'}^{a, v})) ) \\
    b_{r,r'}^{b, v} = \min(ToC - (e_{r^+}^v - g_{r'^- \xrightarrow{} r'^+} + b_{r,r'}^{a, v}), t_{r, r'}^{v, b} (\gamma_s s_{r,r'}^{b, v} + \gamma_f (1-s_{r,r'}^{b, v})) ) \\
    b_{r,r'}^{v} = b_{r,r'}^{a, v} + b_{r,r'}^{b, v}, \; \forall v \in V, \forall r, r' \in \mathcal{R} 
\end{align}

\begin{align}
    e_{{r_0^v}^-}^v = e_{{r_0^v}^+}^v = ToC, \; \forall v \in V, \; \forall r_0^v \in U \\
    e_{r'^+}^v = e_{r'^-}^v - g_{r'^- \xrightarrow{} r'^+}, \; \forall v \in V, \forall r, r' \in \mathcal{R}  \\ 
    e_{r'^-}^v \leq e_{r^+}^v - g_{r^+ \xrightarrow{} r'^-} + b_{r,r'}^v + (1 - y_{r, r'}^v)M, \; \forall v \in V, \forall r, r' \in \mathcal{R}  \\ 
    e_{r'^-}^v \geq e_{r^+}^v - g_{r^+ \xrightarrow{} r'^-} + b_{r,r'}^v - (1 - y_{r, r'}^v)M, \; \forall v \in V, \forall r, r' \in \mathcal{R}  \\ 
    y_{r, r'}^v \cdot SoC_{min} \leq e_{r'^-}^v, \; \forall v \in V, \forall r, r' \in \mathcal{R}  \\ 
    y_{r, r'}^v \cdot SoC_{min} \leq e_{r^+}^v + b_{r,r'}^{a, v}, \; \forall v \in V, \forall r, r' \in \mathcal{R} 
\end{align}

\begin{align}
     s^{a, v}_{r, r'} = s^{b, v}_{r, r'}, \; \forall v \in V, \forall r, r' \in \mathcal{R} \; \; \text{s.t.} \; \; r^+ = r'^- \\
     n_{u} =  \sum_{r \in \mathcal{R}} 1 - \sum_{v \in V} \sum_{r \in \mathcal{R}} y_{r', r}^v  \\
     n_f = \sum_{v \in V} \sum_{r, r' \in \mathcal{R}} (s^{a, v}_{r', r} + s^{b, v}_{r', r}) \cdot 1_{r^+ \neq r'^-} + s^{a, v}_{r', r} \cdot 1_{r^+ = r'^-}
\end{align}

\begin{align*}
     y_{r, r'}^v, s^{a, v}_{r, r'}, s^{b, v}_{r, r'} \in \{0, 1\}\\
     t_{r,r'}^{v, a}, t_{r,r'}^{v, b} \in [0, |T|] \\
     b_{r,r'}^{a, v}, b_{r,r'}^{b, v} \in [0, ToC] \\ 
     e_{r^-}^v, e_{r^+}^v \in [BoC, ToC]
\end{align*}

(12) ensures that each request is served at most by one aircraft, (13) and (14) ensure path coherence (i.e., one aircraft cannot travel multiple routes) and continuity is enforced by (15). 
Charging times between service are constrained by (16)-(18) making sure they do not exceed the maximum time available for charge, accounting for the minimum ground time between successive landing and takeoff when a deadhead is necessary (i.e., $r^+ \neq r'^-$).
(19) and (20) are used to track the electricity bought between two services, respectively after the first service (i.e., $r$) and before the second service (i.e., $r'$). The $\min$ function explicitly ensures the Top of Charge (see Section \ref{sec:batteries}) can never be exceeded. (21) is used to aggregate the electricity bought between $r$ and $r'$, to be included in the objective function (11). (19) and (20) are non-linear constraints due to the minimum operator and the products within it, however both constraints can be linearized using the big-M method at the cost of additional variables. State of Charge is monitored at discrete points : $e_{r^-}$, $e_{r^+}$ which correspond to beginning and ending of flights. The aircraft's State of Charge are tracked using (23)-(25), noting that initial levels (i.e., $\{e_{r^+}, e_{r^-} \mid r \in U\}$) are set to maximum capacity $ToC$ in (22). (26) and (27) ensure that aircraft may not takeoff with a battery charge level below a set threshold, $SoC_{min}$. Charging mode may change only after a landing and therefore must be the same over a connection that does not require a deadhead (28).
Finally, (29) and (30) track the number of fast charges used and unserved demands. Since we allow the charging mode to change only after a flight, the fast charge count differentiates between connections with and without deadheads.


\section{Algorithms}

In this section we introduce new dedicated algorithms to solve the air-taxi problem described in the previous sections. We first present a Beam Search algorithm for the Pooling and Scheduling problem, then a Variable Neighborhood Search for the Routing and Recharging problem. Both algorithms offer some form of incremental solving, allowing to support batched and incremental optimization of sets of demands and requests.

\subsection{Pooling and Scheduling}

The mathematical model presented in Section \ref{sec:pooling} can be solved using a complete solver but can require an important time budget. Besides, as we will see in section \ref{sec:online} we can benefit from reusing past solutions which is typically difficult with a complete solver. We propose a simple and efficient beam search algorithm (UAM-Beam Search) which is able to find good solutions quickly and can solve incremental series of pooling problems by reusing past solutions.

Let us first remark that in any optimal solution to the Pooling and Scheduling problem, constraint (8) and the objective function (3) jointly imply that scheduled departure times $f_k$ are such that $f_k = \max_{i} \{q_{1-\delta}^{\nu_i} \mid x_i^k = 1\}$ for all groups $k$. This means that once group assignments are determined, flight departure times can be deduced without any optimization. Using this, we can characterize solutions using group assignments only (i.e., using variables $x_i^k$), and flight departure times $f_k$ can be deduced. Furthermore, any of these assignments is a partition of the set of demands $\mathcal{D}$, in which each demand is identified by a unique integer identifier \emph{id}, $\{0, \dots, D-1\}$. Exploring the space of partitions of $\{0, \dots, D-1\}$ can in theory be done by building a tree enumerating all possibilities, but could not be done in practice past a certain value of $D$ because of combinatorial explosion. Each level $j$ of such a tree contains all the partitions of the set $\{0, \dots, j\}$ (one partition per node) and is constructed by adding $j$ to the parents partitions of level $j-1$. The root of the tree is initialized with demand $0$ and therefore contains the set $\{\{0\}\}$, level 1 is constructed by including demand 1 in every possible way, $\{\{0, 1\}\}$ or $\{\{0\}, \{1\}\}$, resulting in two nodes. The process is repeated for the next demands until level $D-1$ is reached. Some demands are incompatible with each other due to constraints (7) and (6). These incompatibilities can be evaluated on intermediate nodes, before reaching the last level of the tree. If an intermediate node contains a conflict, it is straightforward that all nodes of the associated sub-tree contain the same conflict and therefore no valid solution is present in the sub-tree. This observation allows us to prune the tree while constructing it to avoid exploring partitions that do not contain valid solutions.

A beam search algorithm consists in extending a search tree under a fixed number of nodes $W$ at each level, instead of extending it under every node, thereby avoiding the exponential explosion of the number of nodes. In fact, this yields a space complexity linear in the depth of the tree, which is in our case equal to the number of demands. At each level, the nodes under which to continue the search are selected according to an \textit{evaluation} function which assigns a value to each node. This evaluation function computes a truncated version of the sum (3) defined in Section \ref{sec:pooling}, considering only the demands that are included in the tree at this stage of the search. We only evaluate nodes which are free of contradiction. Algorithm \ref{alg:beam} presents the pseudo-code of the UAM-Beam Search.

\begin{algorithm}[htbp]
\SetAlgoLined
\textbf{Input:} Beam Width: W, Instance Params: P, Number of demands: D\;
\KwResult{Set of at most W partitions, ranked by quality}
\textbf{Initialize:} Parents $\xleftarrow{}$ \{\{0\}\}, Children $\xleftarrow{}$ \{\}\;
 \For{d=1, \dots, D-1}{
  Children $\xleftarrow{}$ ExtendPrune(Parents, d, P)\;
  Children $\xleftarrow{}$ RetrieveBest(Children, W, P)\;
  Parents $\xleftarrow{}$ Children
 }
 \caption{UAM-Beam Search}\label{alg:beam}
\end{algorithm}

For clarity, instance parameters are noted P and contain all parameters discussed in Section \ref{sec:pooling}: $\alpha_{\mathit{regular}}$, $\alpha_{\mathit{premium}}$, $t_{\mathit{regular}}$, $t_{\mathit{premium}}$, $t^{max}_i$ and lookup tables for arrival time distribution quantiles, expectations and number of passengers, associated with each demand. These are used to determine conflicts and evaluate partition quality. Two auxiliary functions are used. The first one, \texttt{ExtendPrune}, builds the set of all children nodes given a set of parents nodes, a new demand \emph{id} and instance parameters. Children nodes which contains conflicts are immediately pruned. The second function, \texttt{RetrieveBest}, evaluates all nodes produced by \texttt{ExtendPrune} and returns the $W$ best children in ascending order according to the node evaluation function (a truncated version of the sum (3)). The algorithm stops when all the demands were inserted in the tree.

\subsection{Routing and Recharging}
\label{sec:vns-method}
Similarly, the Routing and Recharging problem (Section \ref{sec:routingandcharging}) can be solved with a complete solver. However, as discussed in Section \ref{sec:expe}, this approach does not scale to large instances. Instead, we define a specific Variable Neighborhood Search algorithm, UAM-VNS, to efficiently find good solutions to the Routing and Recharging problem. A general variable neighborhood search contains 3 main phases: shaking, descent and neighborhood change \cite{hansen2017variable}. These phases are repeated until a stopping criterion is met. The shaking part consists in choosing a random neighbor of the incumbent solution, which serves as starting point for the subsequent descent phase. In our case, the descent phase is a Variable Neighborhood Descent (VND). After the descent, a neighborhood change function determines if the incumbent solution is replaced by the one found by VND and which neighborhood to use for the next shaking. The search stops if the time budget has been entirely consumed or if no improvement was made over the last $\eta_{\mathit{improve}}$ iterations. We provide a pseudo-code for our UAM-VNS in Algorithm \ref{alg:vns}, and present details for each part.

\begin{algorithm}[htbp]
\SetAlgoLined
\textbf{Input:} Set of neighborhoods $N^{shake}$ ($\{N_i^{shake}$, $i=1, \dots, n\}$) to be used in the shaking phase, and a set $N^{search}$ ($\{N_j^{search}$, $j=1, \dots, m\}$) to be used in the descent phase, Evaluation function F\;
\KwResult{Solution to the Routing and Recharging Problem}
\textbf{Initialize:} Initial solution $x$, Initial value for penalty parameter $\lambda$\;
 \While{Stopping criterion}{
  $i \xleftarrow{} 1$\;
  \While{$i < n$}{
    \textbf{Shaking:} $ x' \xleftarrow{}$ random neighbor in $N_i^{shake}(x)$\;
    \textbf{Descent:} $x'' \xleftarrow{}$ VND(x', $N^{search}$, $\lambda$, F)\;
    \textbf{Move:} If $F(x'', \lambda) < F(x, \lambda)$, then $x \xleftarrow{} x''$ and $i \xleftarrow{} 1$. Otherwise  $i \xleftarrow{} i + 1$\;
    \textbf{Update:} $\lambda \xleftarrow{}$ UpdatePenalty($\lambda$)
  }
  }
 \caption{UAM-VNS}\label{alg:vns}
\end{algorithm}

\paragraph{Neighborhoods}

We define 4 neighborhood functions, described in Table \ref{tab:neigh}, 3 targeting routes and one targeting charging operations. Neighborhood $N_1$ shifts a request from an aircraft to another aircraft or to no aircraft. When applied to a served request, $N_1$ can either re-assign the request to a different aircraft or to no aircraft, increasing the unserved requests count by 1. When applied to an unserved request, $N_1$ assigns the request to one of the aircraft in the fleet. Neighborhood $N_2$ applies to pairs of requests served by different aircraft and swaps their aircraft. Neighborhood $N_3$ performs a greedy and local optimization of charging mode and charging time immediately before or after a flight, using the golden section search method \cite{10.2307/2032161}. Neighborhood $N_4$ exchanges full routes between aircraft, potentially removing requests that are infeasible due to the initial positions of the aircraft (the fleet of aircraft we are considering is homogeneous but aircraft start at different locations), therefore always returning feasible solutions.

\begin{center}
\begin{table}[htbp]
\begin{tabular}{ l p{12cm} }
 \hline 
 Neighborhood & Description \\
 \hline
  $N_1$ : Shift &  Shift one request $r$ from its current assigned aircraft to another one. \\ 
  $N_2$ : Swap &  Swap two requests that are served by different aircraft. \\ 
  $N_3$ : Recharge & Change one charging duration and mode. \\ 
  $N_4$ : Rotate & Reassign existing paths to different aircraft. \\ 
 \hline
\end{tabular}
\caption{\label{tab:neigh} UAM-VNS Neighborhoods }
\end{table}
\end{center}

\paragraph{Shaking}

The shaking part of the GVNS aims at favoring diversification and helps escaping local minima. It consists in selecting a random point in the neighborhood of the current solution. In our UAM-VNS we alternate between $N_1$ and $N_2$ for this phase, i.e. $N^{shake} = \{N_1, N_2\}$.

\paragraph{Variable Neighborhood Descent}

In the VND phase, we use all 4 neighborhood structures, i.e., $N^{search} = \{N_1, N_2, N_3, N_4\}$, and a best improvement strategy. Neighborhoods are used sequentially to move around the starting point, until a local minimum with respect to all 4 neighborhoods is found. A pseudo-code for this procedure is given in Algorithm \ref{alg:vnd}.

\begin{algorithm}[htbp]
\SetAlgoLined
\textbf{Input:} Set of neighborhoods $N^{search}$ ($\{N_j^{search}$, $j=1, \dots, m\}$) to be used in the descent phase, Evaluation function F, Initial solution $x'$, Initial value for penalty parameter $\lambda$.\;
\KwResult{Local minima with respect to $N^{search}$.}
  j $\xleftarrow{} 1$\;
  \While{$i < m$}{
    \textbf{Local Search:} $x'' \xleftarrow{}$ $argmin_{y \in N_j(x')} F(y, \lambda)$\;
    \textbf{Move:} If $F(x'', \lambda) < F(x', \lambda)$, then $x' \xleftarrow{} x''$ and $j \xleftarrow{} 1$. Otherwise $j \xleftarrow{} j + 1$\;
  }
 \caption{Variable Neighborhood Descent}\label{alg:vnd}
\end{algorithm}

\paragraph{Objective}

During the search we allow for violation of the SoC constraints and penalize for it, yielding the following objective function.

\begin{equation}
    F(x, \lambda) = f(x) + \lambda \cdot h(x)
\end{equation}

Where $f$ is the objective function defined in (11), $x$ is a candidate routing and recharging solution and $h$ is a function that maps any solution to its degree of SoC constraints violation. We define the degree of violation as the sum of the negative slack excess of constraints (26) and (27). The penalty parameter $\lambda$ is dynamic and is updated in the update step of Algorithm 2. We periodically decrease the penalty by a multiplicative constant, $\beta_{decr}$, when the last few visited solutions were all valid and increase it by and additive constant, $\beta_{incr}$, otherwise \cite{DBLP:conf/ijcai/HamadiJS09, hamadi2008manysat}.

\section{Results and Discussions}\label{sec:expe}

\subsection{General Settings}

Several air-taxi services are under active development
\cite{Volocopter2019,ehang2020,whitepapercommunityelevate2020,UberWP2016}.
We used insights from these studies to design our experiments. 
For the Pooling and Scheduling problem, synthetic demand scenarios are generated from a passenger commuting model and used to evaluate our methods. For the Routing and Recharging problem, we design 18 infrastructures characterized by their number of vertiports and aircraft, and for each of them, 3 synthetic request scenarios are generated.
All experiments were run on a MacBook Pro laptop equipped with a 6 cores 2,6 GHz Intel Core i7 processor with 32 Gb of RAM,  running MacOS catalina v10.15.7. In all Gurobi v9.1.0 experiments, we set an optimality gap of $10^{-3}\%$ and up to 4 threads, the rest of the parameters were left to their default values. Each algorithm had a 30 minutes time-out for each run.

\subsection{Evaluation, Pooling and Scheduling}

\subsubsection{Specific settings}\label{sec:settings-pooling}

To evaluate our pooling and scheduling method, we generate synthetic demands. Each demand concerns a specific route (i.e., a pair of origin and destination vertiports) and is characterized by a number of passengers, an arrival time distribution, a constraint on arrival time and a class (Regular or Premium), see Section \ref{sec:pooling}. In this section, we only consider one route, therefore we only focus on demand attributes described in Table \ref{tab:notations-demand}. Throughout this section, unless specified otherwise, Pooling and Scheduling instances are generated by drawing each demand $i$ from the following stochastic model:

\begin{itemize}
    \item The number of passengers $w_i$ is sampled from the distribution $\{1:70\%, 2:20\%, 3:5\%, 4:5\%\}$;
    \item As described in Section \ref{sec:pooling}, our pooling and scheduling model only uses the mean and quantile $q_{1-\delta}^{\nu_i}$ of a passenger's arrival time distribution $\nu_i$. The mean, expressed in minutes, is sampled from a mixture between two normal distributions, centered around 8:30 am and 5:00 pm with a standard deviation of 20 minutes, and a uniform distribution over the full service period. The weights of the mixture are respectively $\frac{1}{2}, \frac{1}{3} \; \text{and} \; \frac{1}{6}$. This approximates a commuting pattern. The $1-\delta$ quantile is modelled as an offset expressed in minutes relative to the mean arrival time, and sampled from the distribution $\{3:40\%, 5:50\%, 7:10\%\}$, describing 3 degrees of confidence on the arrival time at the source vertiport.
    \item A constraint on arrival time at destination vertiport is present with probability $20\%$. When present, it is modeled as a maximum offset $t^{\mathit{max}}_i$ expressed in minutes relative to $q_{1-\delta}^{\nu_i}$, sampled uniformly from $\{10, 15, 20\}$.
    \item The class $c_i$ is Premium with probability $20\%$ and Regular with probability $80\%$.
\end{itemize}

Additionally, and unless specified otherwise, the maximum allowed expected waiting times for regular demands (i.e., $t_{\mathit{regular}}$) is set to 25 minutes and to 15 minutes for premium demands (i.e., $t_{\mathit{premium}}$). We use weights $\alpha_{\mathit{regular}} = 1$ and $\alpha_{\mathit{premium}} = 2$ (see objective function (3)) to penalize for the expected waiting time of both classes. Aircraft capacity is set to 4 passengers.\\

Our UAM Beam Search is configured with a beam-width parameter of $1000$ for all experiments reported in this section. The performance analysis is conducted for several number of demands $D$. For each value of $D$, $20$ instances are generated and solved to account for variability in our generation process (see previous paragraph) and we report means ($\pm$ standard deviations) for several metrics. These are the number of requests (\#Requests) created after pooling and scheduling, the number of demands per created requests (Requests Load) and expected waiting times for premium and regular demands (WT Premium, WT Regular).
\newline

\begin{table}[htbp]
\centering
\resizebox{\textwidth}{!}{
\renewcommand{\arraystretch}{1.7} 
 \begin{tabular}{clccccccc}
 \toprule
 & & \multicolumn{7}{c}{Number of Demands D}
 \\\cmidrule(r){3-9}
& &  20  & 25 & 30 & 35  & 40 & 45 & 50  \\
\hline
\parbox[t]{2mm}{\multirow{5}{*}{\rotatebox[origin=c]{90}{\textbf{Gurobi 9.1.0}}}} 
& \#Requests & 10.1 $\pm$ 0.9 & 12.8 $\pm$ 1.2 & 14.0 $\pm$ 1.0 & 16.0 $\pm$ 1.1 & 17.8 $\pm$ 1.0 & 19.2 $\pm$ 1.7 & 21.4 $\pm$ 1.6 \\
& Requests Load  & 2.6 $\pm$ 1.1 & 2.6 $\pm$ 1.2 & 2.9 $\pm$ 1.2 & 2.9 $\pm$ 1.2 & 3.0 $\pm$ 1.2 & 3.1 $\pm$ 1.2 & 3.1 $\pm$ 1.1 \\
& WT Premium (') & 6.4 $\pm$ 3.4  & 6.1 $\pm$ 3.0  & 6.1 $\pm$ 3.2 & 6.0 $\pm$ 2.8 & 6.2 $\pm$ 3.0 & 6.1 $\pm$ 2.9 & 6.0 $\pm$ 2.7 \\
& WT Regular (') & 8.2 $\pm$ 5.7 & 8.1 $\pm$ 5.7 & 8.3 $\pm$ 5.4 & 8.2 $\pm$ 5.5 & 7.6 $\pm$ 4.8 & 8.0 $\pm$ 5.2 & 7.7 $\pm$ 4.8 \\
& Time ('') & 0.3 $\pm$ 0.1 & 1.0 $\pm$ 0.3 & 1.8 $\pm$ 0.4 & 3.4 $\pm$ 0.9 & 6.7 $\pm$ 2.0 & 22.0 $\pm$ 19.6 & 104.3 $\pm$ 142.8 \\
\hline 
\parbox[t]{2mm}{\multirow{5}{*}{\rotatebox[origin=c]{90}{\textbf{UAM-Beam Search}}}} 
& \#Requests & 10.1 $\pm$ 0.9 & 12.8 $\pm$ 1.2 & 14.0 $\pm$ 0.9 & 16.2 $\pm$ 1.1 & 18.1 $\pm$ 1.1 & 19.7 $\pm$ 1.8 & 21.9 $\pm$ 1.6\\
& Requests Load  & 2.6 $\pm$ 1.2 & 2.6 $\pm$ 1.2 & 2.8 $\pm$ 1.2 & 2.9 $\pm$ 1.2 & 2.9 $\pm$ 1.2 & 3.0 $\pm$ 1.2 & 3.1 $\pm$ 1.2\\
& WT Premium (') & 6.4 $\pm$ 3.4 & 6.2 $\pm$ 3.0 & 6.0 $\pm$ 3.0 & 6.1 $\pm$ 2.8 & 6.4 $\pm$ 3.2 & 6.3 $\pm$ 3.0 & 6.2 $\pm$ 2.8\\
& WT Regular (') & 8.2 $\pm$ 5.7 & 8.1 $\pm$ 5.7 & 8.3 $\pm$ 5.5  & 8.4 $\pm$ 5.8 & 8.2 $\pm$ 5.5 & 8.5 $\pm$ 5.5 & 8.3 $\pm$ 5.5\\
& Time ('') & 0.02 $\pm$ 0.00 & 0.04 $\pm$ 0.00 & 0.06 $\pm$ 0.00 & 0.09 $\pm$ 0.00 & 0.12 $\pm$ 0.00 & 0.17 $\pm$ 0.00 & 0.23 $\pm$ 0.02\\
 \bottomrule
\end{tabular}}

\caption{Comparison between Gurobi 9.1.0 and our Beam Search for Pooling and Scheduling.
}
    \label{tab:poolingbench}
\end{table}

\subsubsection{Performance analysis}

Table \ref{tab:poolingbench} reports results obtained when comparing the performance of Gurobi 9.1.0 to that of our UAM Beam Search when ran over a set of simulated instances of the pooling and scheduling problem. UAM Beam Search consistently finds the optimal solution for small sizes (i.e., $D = 20$, and $25$) and finds solutions close to the optimal for all other sizes while being orders of magnitude faster than Gurobi. Furthermore, as instances grow in size, Gurobi shows a high variance in the resolution time while our UAM Beam Search remained stable. 

On average, pooling and scheduling results in the creation of one request for every two demands. As $D$ increases and more demands are simulated, the average requests load increases since more pooling opportunities naturally appear. This increase also leads to smaller average and variance in expected waiting times for both passenger classes, although this is noticeable only in the optimal solutions obtained by Gurobi. 

As intended, premium demands wait less than regular ones on average. We also notice that, as a by-product of the optimization, premium demands enjoy a lower variability (standard-deviation) in their expected waiting time.

\begin{table}[htbp]
\centering
 \begin{tabular}{lcc}
 \toprule
 & Peak hours (7:30 - 9:30 am and 4:00 - 6:00 pm) & Off-Peak hours   \\
\hline

WT Premium (') & 6.6 $\pm$ 3.3 & 5.4 $\pm$ 2.8 \\
WT Regular (') & 8.7 $\pm$ 5.3 & 7.4 $\pm$ 5.7 \\
\bottomrule
\end{tabular}
\caption{Quality of Service during peak and off-peak hours.}
\label{tab:peakhours}
\end{table}

We used synthetic data derived from a commuting pattern which creates peak hours and off-peak hours by design. Consequently, we could expect different QoS levels during these different periods. To evaluate this, we generate 100 demands and we measure the expected waiting time during peak and off-peak hours for both classes, after pooling and scheduling. The experiment is repeated 20 times. The results, presented in Table \ref{tab:peakhours}, show that the QoS in on average better during off-peak hours than during peak hours for both classes while still preserving a better QoS for premiums.

\subsubsection{Passenger-centric analysis}

The results of Table \ref{tab:poolingbench} show a difference in waiting times between regulars and premiums, the latter being favored. However it is not clear how the QoS varies between same-class groups and mixed-class groups. To investigate these points in more details and gain more insights, we conduct several experiments.

\paragraph{Group composition and expected waiting time}
Mixing demands with different classes into the same groups might have an effect on their respective waiting times, due to their different constraints and weights. To investigate this point, we evaluate the effect of falling into a mixed group (i.e., containing both classes) from the point of view of each class, considering different values for $\alpha_{\mathit{regular}}$ with respect to a fixed $\alpha_{\mathit{premium}}$. To do so, we generate demands according to the demand model described in introduction of this section and fix $\alpha_{\mathit{premium}} = 1$. We measure the expected waiting times for four categories of demands: regular demands pooled with regular demands only, premium demands pooled with premium demands only, regular demands pooled with both regular and premium demands, and premium demands pooled with both regular and premium demands. Demands that are alone in their group naturally have a minimal expected waiting time. To avoid bias in this experiment, we consider only groups which contain at least 2 demands. For each value of $\alpha_{\mathit{regular}}$ ranging  from $0$ (i.e., no importance is given to minimizing regular expected waiting time) to $1$ (i.e., minimizing both classes' expected waiting time is equally important) with a step of $0.1$, we generate 100 demands and run our UAM-Beam Search. For each value of $\alpha_{\mathit{regular}}$, the experiment is repeated 20 times to account for the randomness of our demand model and we obtain the distributions presented in Figure \ref{fig:params-study-mixity}. We repeat this experiment for several proportions of premium demands (20\%, 50\% and 80\%) representing scenarios in which they are a minority, equally frequent as regular demands or a majority. Observing the results in Figure \ref{fig:params-study-mixity}, we can make the following remarks:

\begin{figure}[htbp]
\begin{subfigure}{1\textwidth}
  \centering
  \includegraphics[width=1\linewidth]{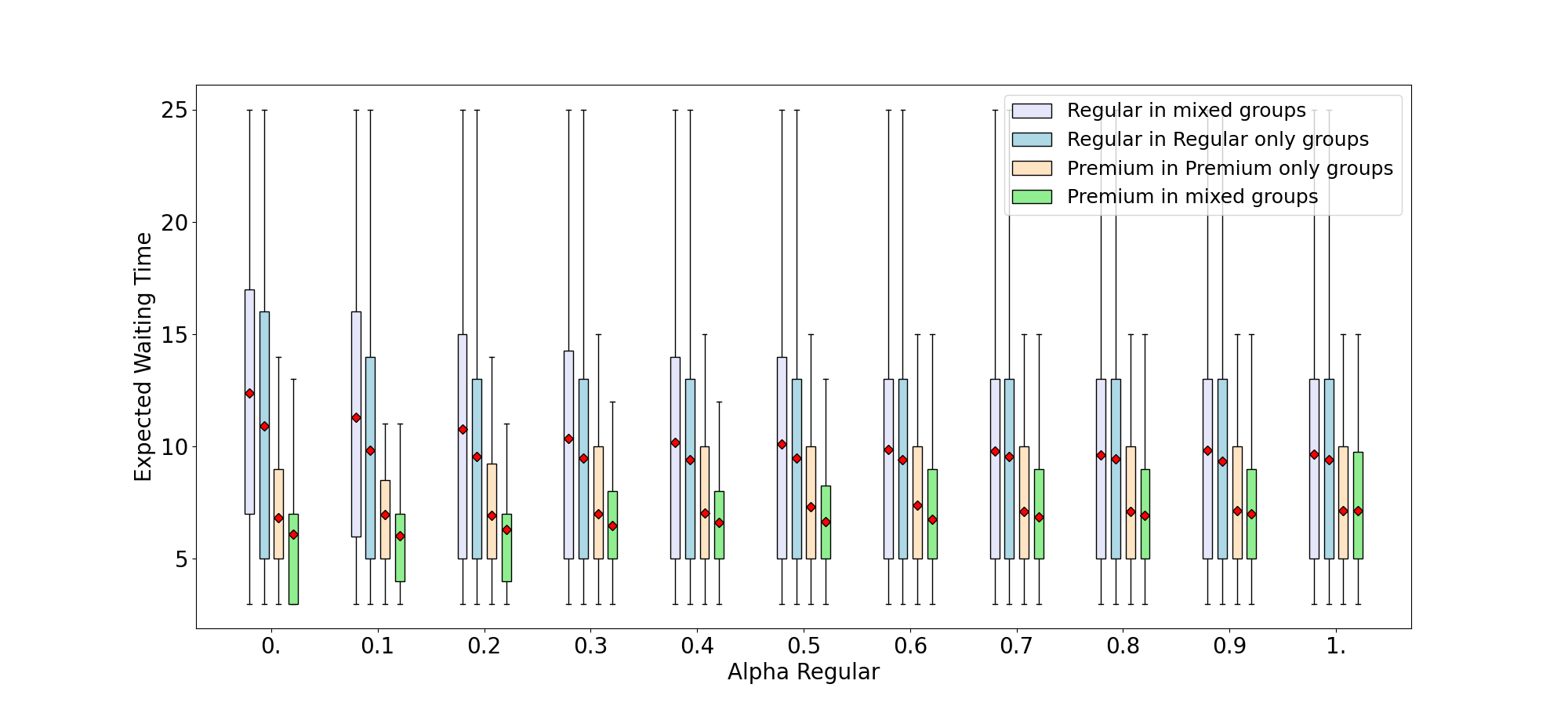}
  \caption{20\% Premium}
  \label{fig:sfig1}
\end{subfigure}\\
\begin{subfigure}{1\textwidth}
  \centering
  \includegraphics[width=1\linewidth]{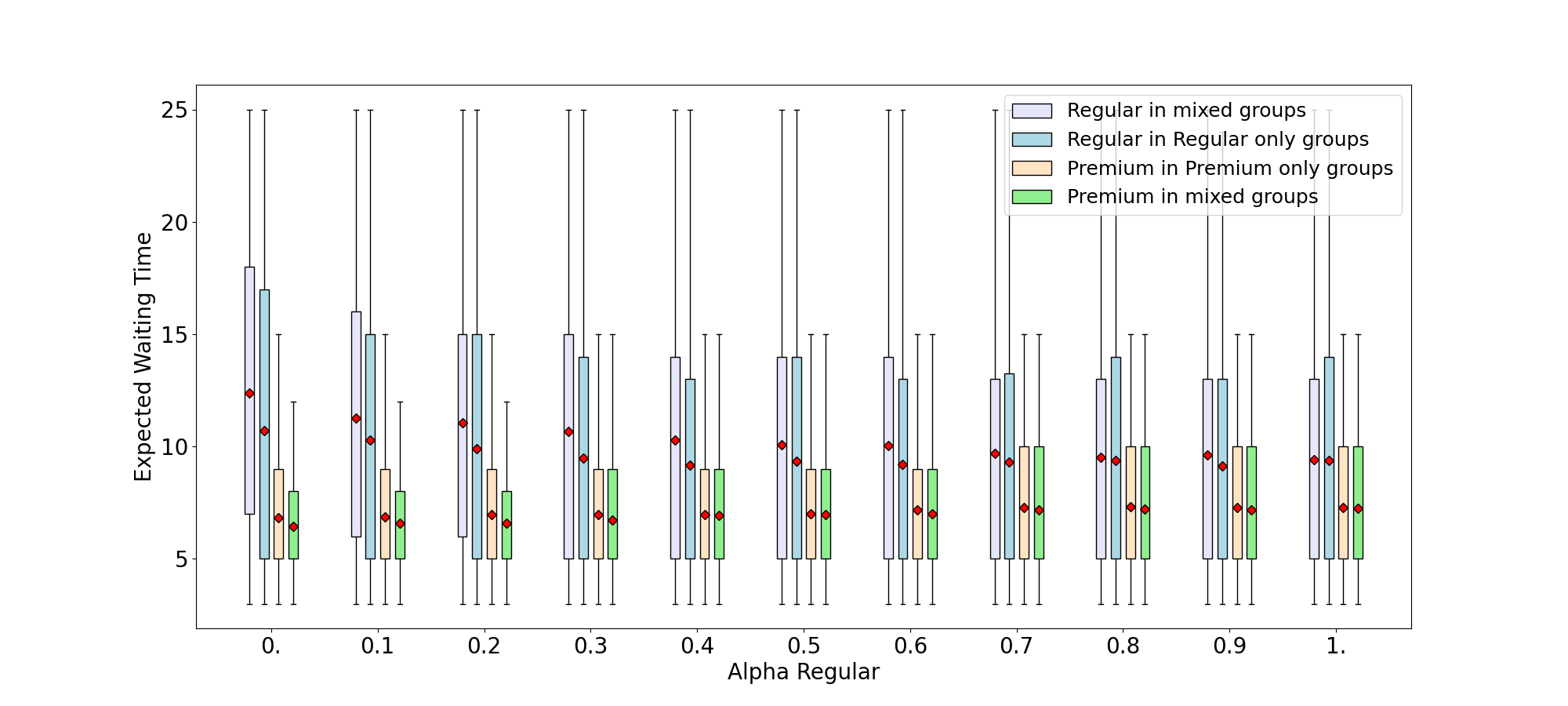}
  \caption{50\% Premium}
  \label{fig:sfig2}
\end{subfigure}\\
\begin{subfigure}{1\textwidth}
  \centering
  \includegraphics[width=1\linewidth]{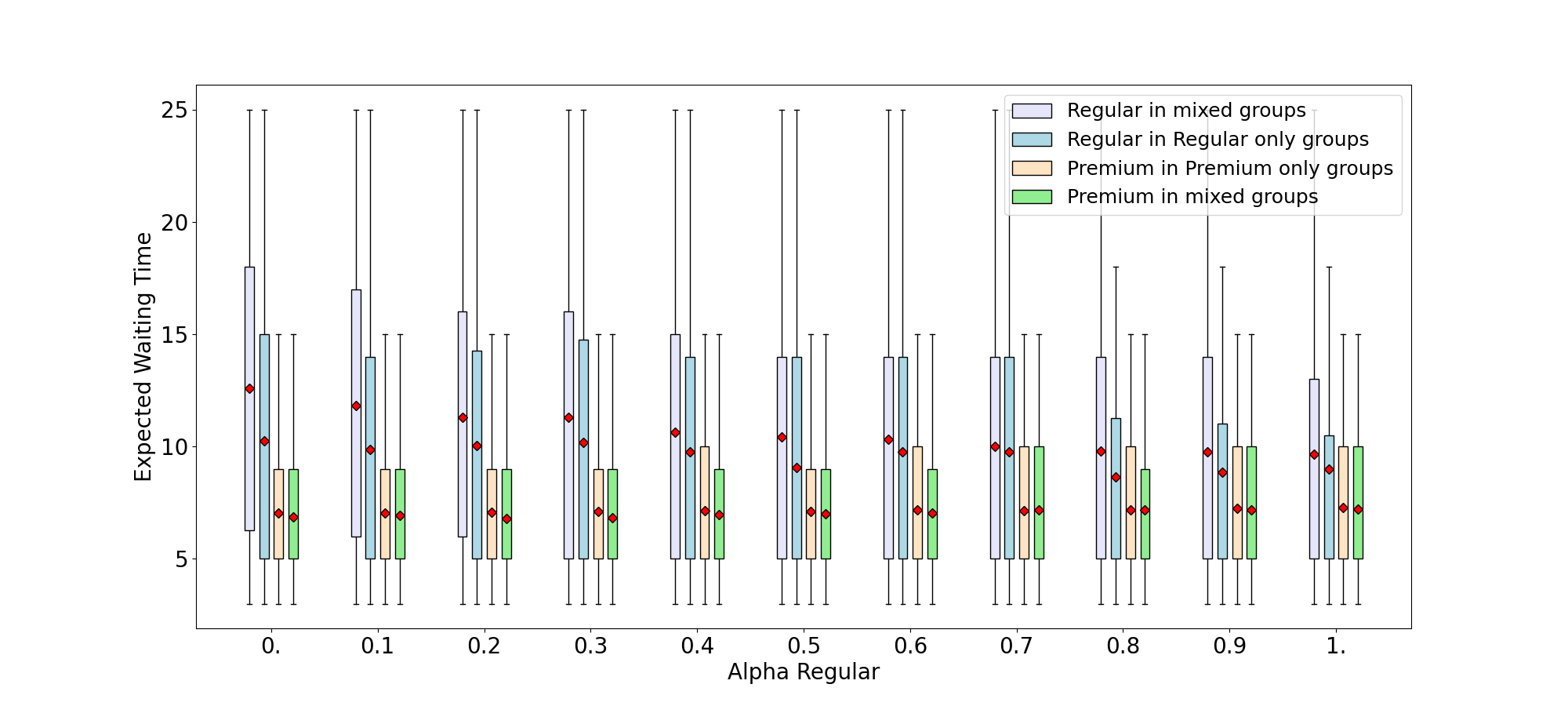}
  \caption{80\% Premium}
  \label{fig:sfig3}
\end{subfigure}
\caption{Distributions of expected waiting times for regular and premium demands clustered by their group composition. Average values are marked in red.
The experiment is repeated with different proportions of premium demands: 20\%, 50\% and 80\%.}
\label{fig:params-study-mixity}
\end{figure}

\begin{itemize}
    \item Overall we observe lower expected waiting times and lower variability for premium than for regular demands, even when $\alpha_{\mathit{regular}} = 1$, i.e., when regular waiting time has the highest importance.
    \item Regular demands do not benefit from being pooled with premiums. On the contrary, they experience longer waiting times. This is particularly the case when $\alpha_{\mathit{regular}} = 0$, and even if the effect wears off as $\alpha_{\mathit{regular}}$ increases (i.e., when giving more importance to the regulars' waiting time), it persists for $\alpha_{\mathit{regular}} = 1$ when the proportion of premium is high (80\%).
    \item Premium demands on the other hand benefit from being pooled with regulars and have better waiting times than when pooled with other premiums. This is particularly marked for a low proportion of premiums (20\%) and low values of $\alpha_{\mathit{regular}}$. This effect wears off when premiums become prominent in the population, even for low $\alpha_{\mathit{regular}}$. The same is observed when $\alpha_{\mathit{regular}}$ is increased, regardless the proportion of premiums.
    \item When $\alpha_{\mathit{regular}}$ is small, we observe in addition that there is a more notable difference in QoS between regular users pooled with premium users versus regular users pooled with only regular users. The same phenomenon is observed on premium demands, albeit with less magnitude.
\end{itemize}

The first remark comes in confirmation of what was already observed in Table \ref{tab:poolingbench}, the persistence of the difference as $\alpha_{\mathit{regular}}$ increases being explained by the tighter bound $t_{\mathit{premium}}$ on the expected waiting time for premium demands. The second and third remarks, which are less obvious, suggest that premiums enjoy a better quality of service (QoS) at the expense of regulars. The fourth remark shows fairness issues within the regular and premium classes when $\alpha_{\mathit{regular}}$ is small. We propose the following explanation for the last three observations:

\begin{itemize}
    \item[] The maximum expected waiting time for premium demands is constrained (see equation (7)) with a tighter bound $t_{\mathit{premium}}$. Therefore, it is less likely for a premium demand to be pooled with another demand which \textit{arrives later} than it is for a regular demand. As consequence, premium demands are more likely to be the \textit{last arriving} in their group, thus enjoying lower waiting times. On the other hand, when regular demands are pooled with premiums, they are less likely to be the \textit{last arriving} and more likely to wait more. Since regular demands are more flexible, premium demands can even benefit from being pooled with them rather than with another premium which also has tight constraints. This benefit disappears as the proportion of regular demands decreases.
\end{itemize}

\begin{figure}[htbp]
\begin{subfigure}{\textwidth}
  \centering
  \includegraphics[width=1\linewidth]{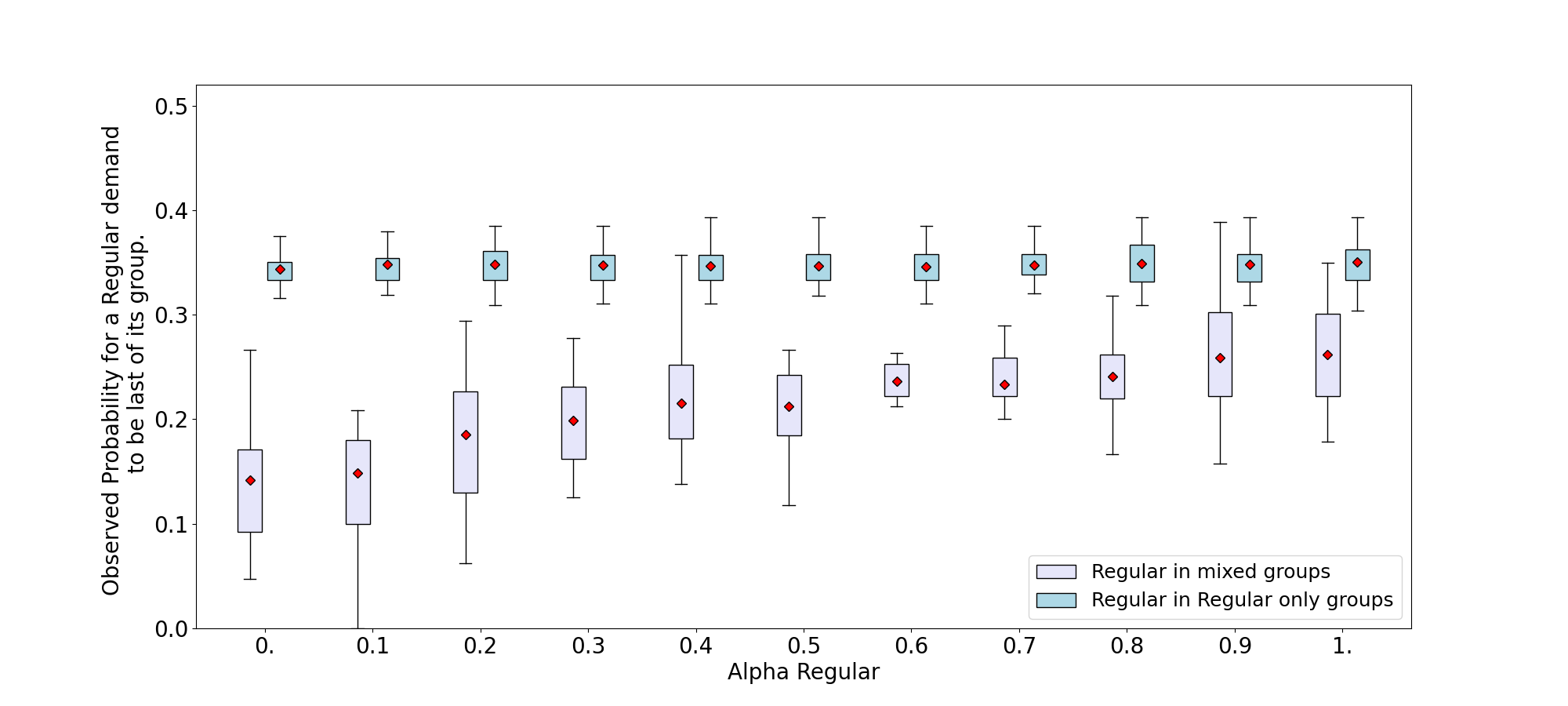}
  \caption{20\% Premium}
  \label{fig:last-arrival-20}
\end{subfigure}\\
\begin{subfigure}{\textwidth}
  \centering
  \includegraphics[width=1\linewidth]{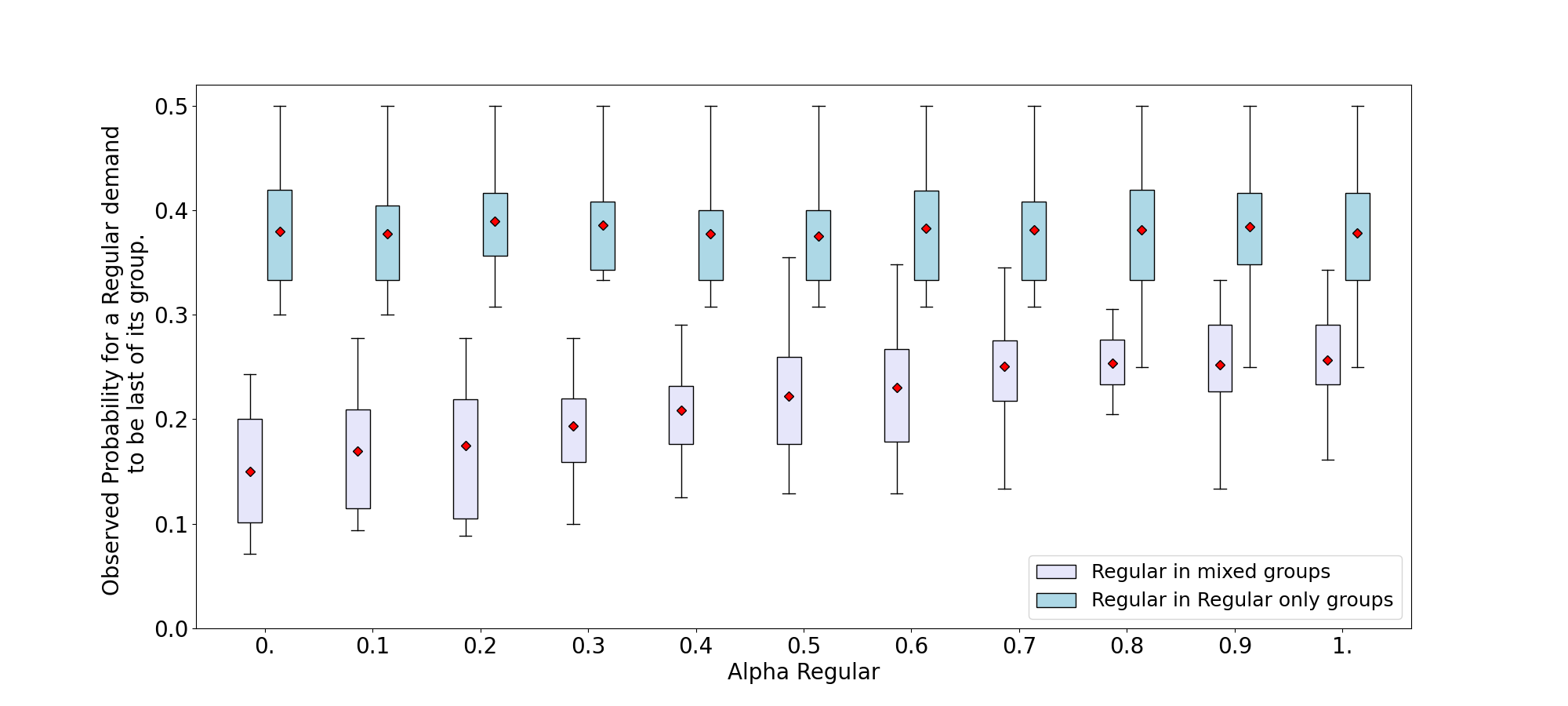}
  \caption{50\% Premium}
  \label{fig:last-arrival-50}
\end{subfigure}\\
\begin{subfigure}{\textwidth}
  \centering
  \includegraphics[width=1\linewidth]{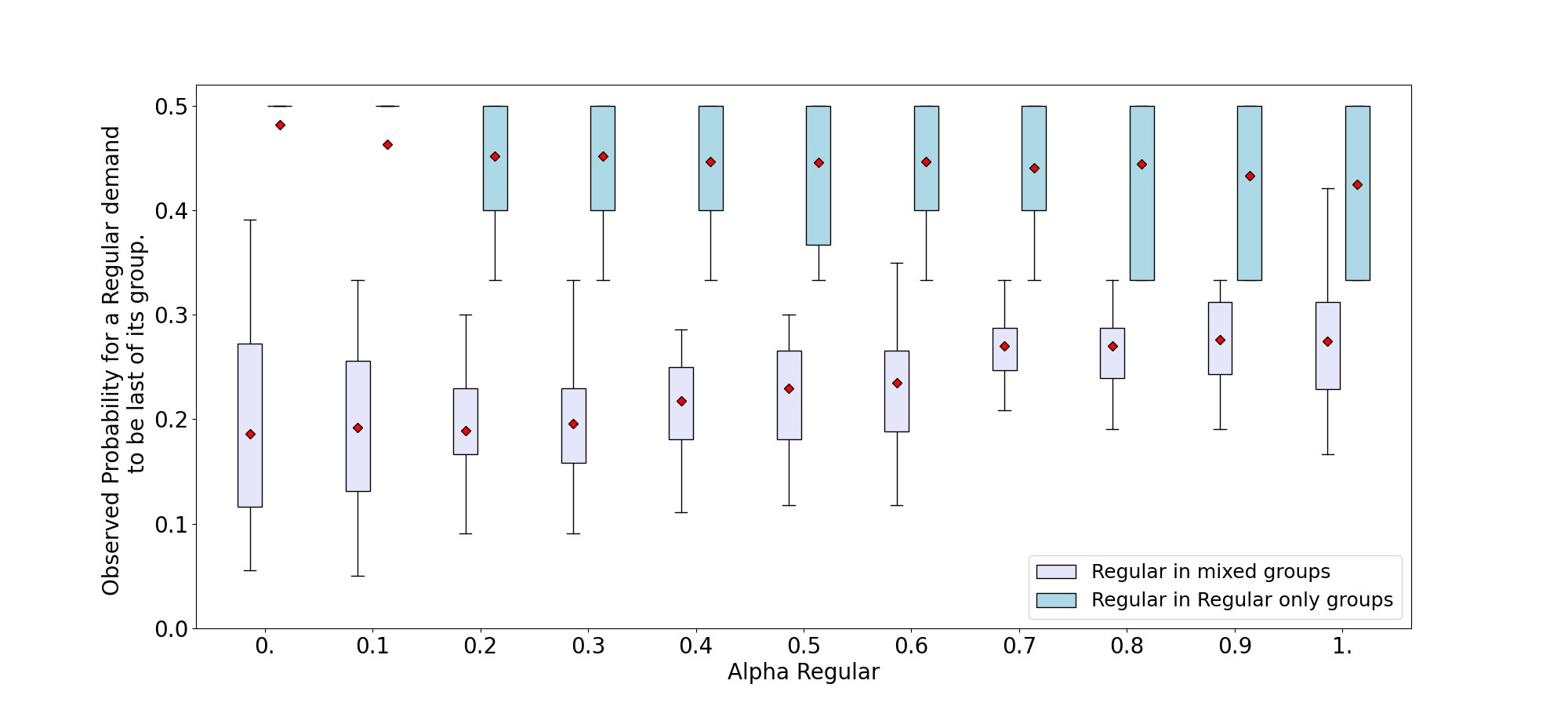}
  \caption{80\% Premium}
  \label{fig:last-arrival-80}
\end{subfigure}
\caption{Distributions of observed probabilities for a regular demand to be the last arriving in its group, clustering by presence of premium demands within the group. 
For each weight $\alpha_{\mathit{regular}}$, 20 samples of size 100 demands were generated to obtain these distributions.
}
\label{fig:params-study-mixity-last-arrival}

\end{figure}

\paragraph{Likelihood of arriving last in a group depending on class}
To verify our previous hypothesis, we design the following experiment. We use the same data generation process as in the previous paragraph. For each value of $\alpha_{\mathit{regular}}$, 100 demands are generated, pooled and scheduled using our UAM Beam Search algorithm. In the obtained solution, we measure the probability (as a frequency) for any regular demand to be the last arriving in its group. As previously, this experiment was repeated 20 times, and we obtain the distributions presented in  Figure \ref{fig:params-study-mixity-last-arrival}. These results show that regular demands are less likely to be \textit{arriving last} in their group when the group contains at least one premium demand. The gap between \emph{regular pooled with premiums} and \emph{regular pooled with regulars} decreases as $\alpha_{\mathit{regular}}$ increases towards $1$. Again, this gap was not completely closed when $\alpha_{\mathit{regular}} = 1$ because premiums still benefit from a tighter bound $t_{\mathit{premium}}$. 

\paragraph{Impact of classes relative differences on QoS}
These two experiments show that in this premium/regular system, regular demands do not benefit from sharing a flight with premiums on average, on the contrary, they experience longer waiting times when this happens.

\begin{figure}[htbp]
\begin{subfigure}{0.5\textwidth}
  \centering
  \includegraphics[width=\linewidth,trim={8cm 2cm 6cm 4.6cm},clip]{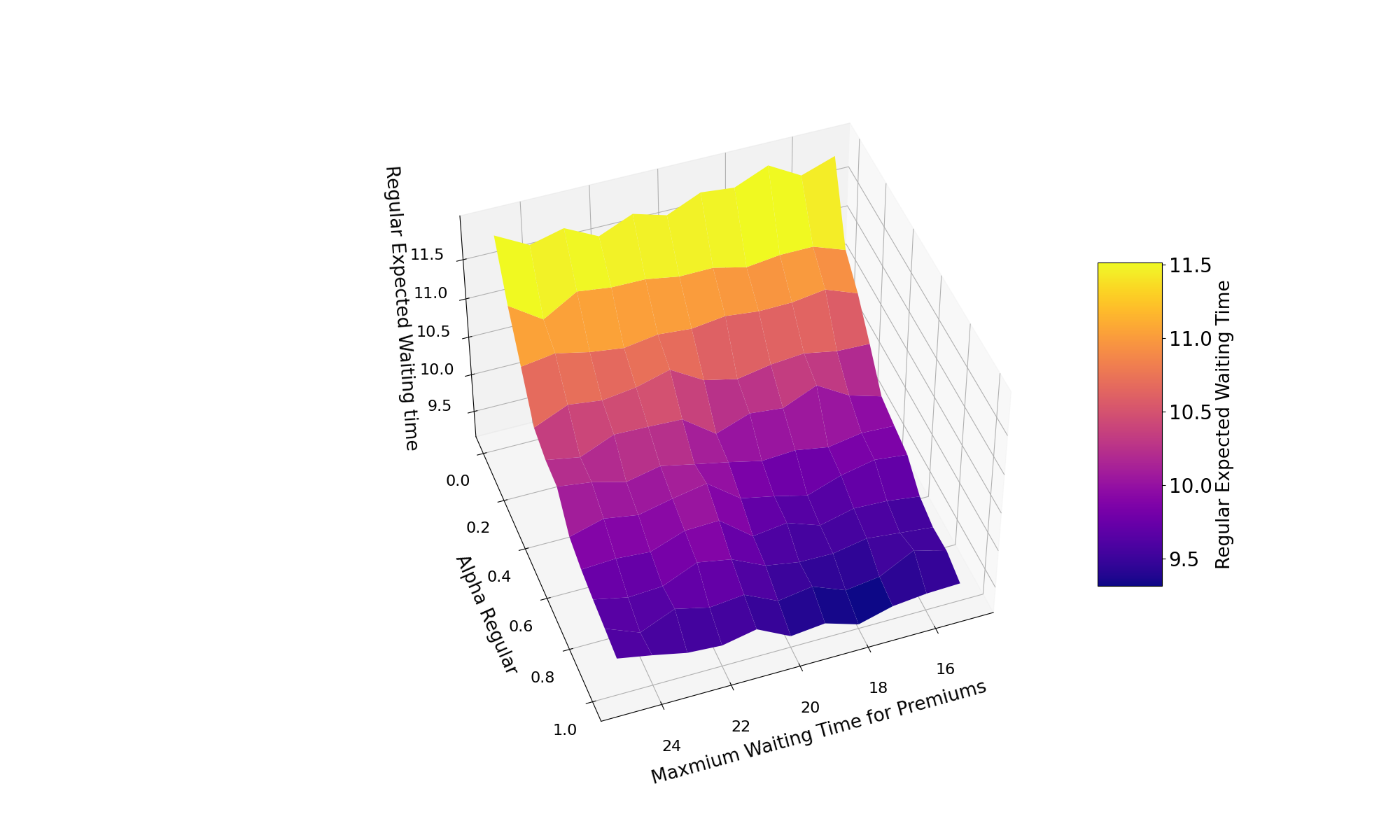}
  \caption{Regular demands}
  \label{fig:hard-constraints-premium-onreg}
\end{subfigure}%
\begin{subfigure}{0.5\textwidth}
  \centering
  \includegraphics[width=\linewidth,trim={8cm 2cm 6cm 4.6cm},clip]{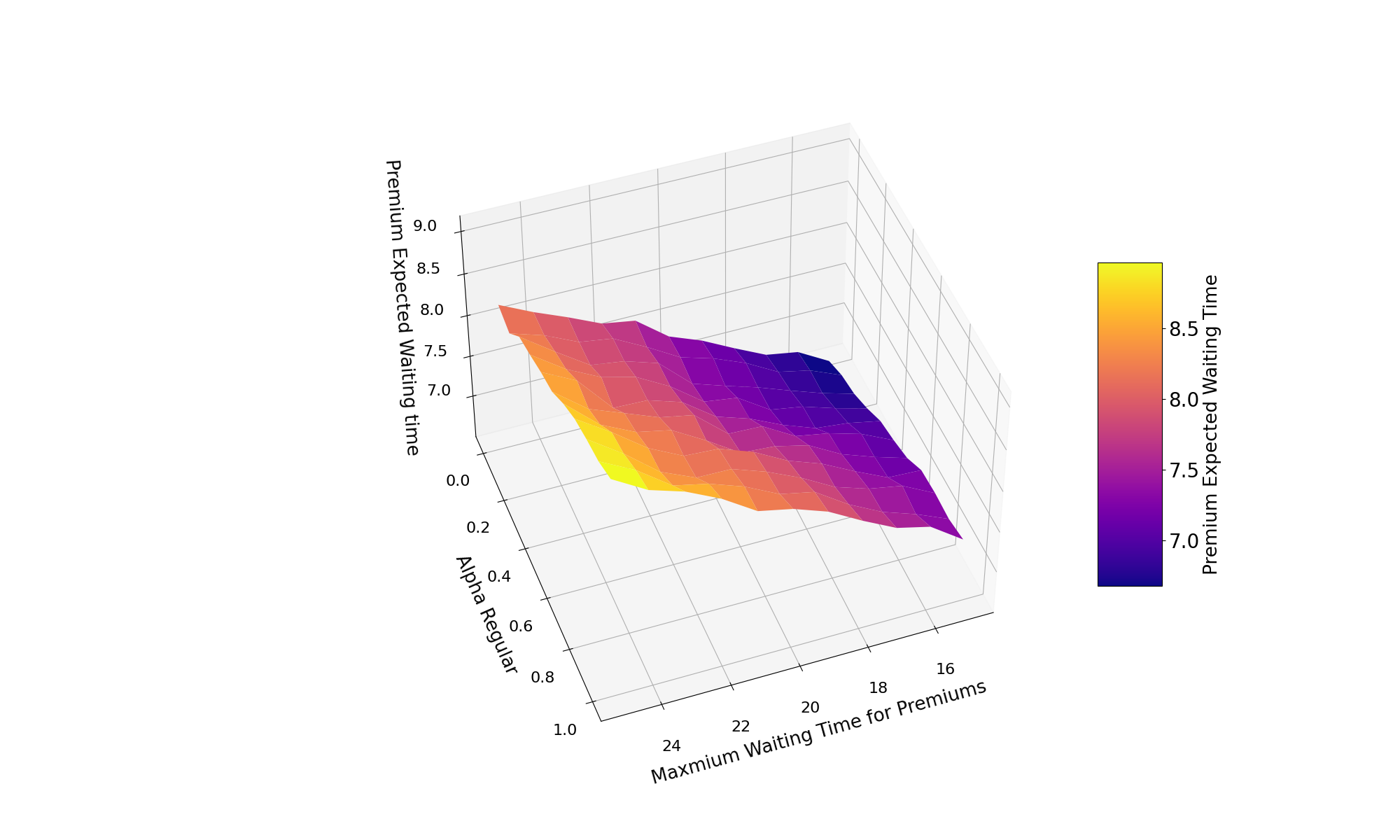}
  \caption{Premium demands}
  \label{fig:hard-constraints-premium-onpre}
\end{subfigure}
\caption{Average expected waiting time for regular (top) and premium (bottom) demands as $t_{\mathit{premium}}$ and $\alpha_{\mathit{regular}}$ both vary. Similar to previous experiments, $\alpha_{\mathit{premium}} = 1$ and we generate $D=100$ demands 20 times for each combination of parameters, and average values of expected waiting times are reported. A proportion of 50\% premiums is used.}
\label{fig:params-study-maxwp}

\end{figure}

However, it is still possible that regular demands experience a better quality of service when premium are present in the system than when regulars and premiums are not differentiated. To test this hypothesis, we use the same generation process as in the previous paragraphs. We measured the average expected waiting times for both classes as $t_{\mathit{premium}}$ and $\alpha_{\mathit{regular}}$ vary respectively between 15 and 25 minutes, and between 0 and 1. We keep $t_{\mathit{regular}} = 25$ and $\alpha_{\mathit{premium}} = 1$. At one extreme ($t_{\mathit{premium}} = 25$ and $\alpha_{\mathit{regular}} = 1$) we have no distinction between the two classes and at the other extreme ($t_{\mathit{premium}} = 15$ and $\alpha_{\mathit{regular}} = 0$) the distinction is maximal. Results of this experiment are reported in Figure \ref{fig:params-study-maxwp}.

When no distinction is made, Figure \ref{fig:params-study-maxwp} confirms that both classes experience a similar average expected waiting time of around 9 minutes. As $t_{\mathit{premium}}$ decreases, \textit{only} the premium class benefits and sees its expected waiting time decrease, regardless of what the value of $\alpha_{\mathit{regular}}$ is. On the other hand, we see no benefit for the regular class which expected waiting time is \textit{solely} influenced by $\alpha_{\mathit{regular}}$. Figure \ref{fig:params-study-maxwp} is obtained using a proportion of premium demands of 50\%, but we observed similar results for 20\% and 80\%.

\subsubsection{Discussion} 

In this section, we evaluated the quality of our UAM-Beam Search by comparison with optimal solutions found by Gurobi, and showed near-optimal solutions can be found quickly on large pooling and scheduling problems. We found that premium demands enjoy a lower average and lower variability in expected waiting times than regular demands. While the lower variability for premium was not an explicit optimization objective, it could be monetized through appropriate premium pricing since it has been shown that users value reliability in shared urban transportation \cite{alonso2020value}.

We provided a first insight into the effect of some model parameters and found that regular demands do not benefit from the presence of premiums in the system. Furthermore, regulars might experience reduced QoS when pooled with premiums as opposed to being pooled with other regulars, which can raise fairness issues within the regular class itself. Our experiments also show under which conditions these issues can arise, and we proposed and verified an explanation for it (the same fairness remark exists for premium demands as well, albeit less intense). 

The pooling and scheduling model and these experiments provide valuable insight into the design-space of an air-taxi service by characterizing the influence of classes and their relative QoS variations as a function of model parameters. Service designers can tune these parameters to match the QoS to the service they want to provide.

\subsection{Evaluation, Routing and Recharging}
\label{sec:vns-eval}
\subsubsection{Specific settings}

To evaluate our UAM-VNS for the routing and recharging problem, we generate several problem instances. A problem instance is characterized by a tuple (Aircraft, Vertiports, Requests) which gives the number of aircraft, vertiports and requests in the instance. Table \ref{tab:params} describes the parameters associated with vertiports, aircraft, recharging actions and requests. For each vertiport, parameters are randomly generated using a uniform distribution over the set of possible values described in the table. Requests are generated uniformly over the day and the routes.

\begin{itemize}
    \item Vertiports are characterized by their landing fees (in \$). These costs are extrapolated from actual heliport operational charges, and expected UAM infrastructure costs. Vertiports are connected through routes - see Section \ref{sec:routes} - with associated flying time and  operating hours.

    \item Aircraft are characterized by their operating cost. This value covers pilot and maintenance costs, insurance costs and taxes, and is given per minute. Battery packs have an associated minimum and maximum state of charge - see Section \ref{sec:batteries}. Minimum ground time represents the time required for passengers ingress, egress and pre-flight checks.

    \item Recharging operations have two modes, slow and fast, associated with different charging rates - see Section \ref{sec:batteries}.
    
    \item For requests values, we use service values proportional to the number of passengers.

\end{itemize}

\begin{table}[htbp]
\centering
\begin{tabular}{ l c } 
 \hline
 \textbf{Vertiports} & \\
 Landing Fees (\$) & 30, 40, 80 \\ 
 Fly time (minute) & 10, 15, 20  \\
 Operating times (hours) & 7am - 7pm\\
 \hline
 \textbf{Aircraft} & \\
 Operating cost $\eta$ (\$/minute) & 34  \\
 ToC & 92 \\
 $SoC_{min}$ & 55 \\
 Minimum ground time $\Delta$ (minute) & 10\\
 \hline
 \textbf{Recharging} & \\
 Slow charging rate (SoC units / minute) & 1  \\
 Fast charging rate (SoC units / minute) & 2 \\

 \hline
\end{tabular}
\caption{Routing and Recharging instance parameters.}
    \label{tab:params}
\end{table}

For each pair (Aircraft, Vertiport), we use three demand scenarios. A \textit{low} demand one with 5 requests per aircraft on average, an \textit{intermediate} scenario with 10, and a \textit{high} demand case with 15. With the previous settings, we produce 18 tuples which characterize different aircraft, vertiports and requests configurations. Tuples are tagged with their respective demand scenario class (l, i, h). In the UAM-VNS algorithm, the penalty parameter $\lambda$ (see Section \ref{sec:vns-method}) is updated every 20 iterations using $\beta_{incr} = 100$ and $\beta_{decr} = 4$ with an initial value of $0$. We give a time budget of 30 minutes for each method. Gurobi can stop early if an optimal solution is found. UAM-VNS can stop early if a stopping criterion is met (search stagnation threshold).

We compare the performance of our UAM-VNS to that of Gurobi 9.1.0. For each method, we measure the percentage of served requests (\% Service), the number of recharging operations (\#Charge) and the proportion of fast charges (\% Fast). Additionally, we measure the average cost per served request (Cost per Service, CPS) for each method. To compare both methods' CPS, we compute the relative difference between the CPS obtained using UAM-VNS to that obtained using Gurobi, negative values indicating a lower (better) CPS for the UAM-VNS.  

To account for variability in the UAM-VNS, it is run 20 times for each instance with different starting points and mean values ($\pm$ standard deviation) are reported for each metric in Table \ref{tab:benchmark}. In this table, time out is indicated by T/O and memory out by M/O. For Gurobi optimality is not guaranteed when the time out is reached and only the last best solutions found before time out are reported.

\begin{table}[htbp]
\centering
\resizebox{\textwidth}{!}{
 \begin{tabular}{lccccccc}
\toprule
 &\multicolumn{3}{c}{\textbf{Gurobi 9.1.0}}
&
\multicolumn{3}{c}{\textbf{UAM-VNS}}
\\\cmidrule(r){2-4}\cmidrule(l){5-7}
(A, V, R) & \% Service & \#Charge (\% Fast) & Time ('') & \% Service & \#Charge (\% Fast) & Time (s) & CPS Diff \\
\hline
(3, 3, 15)l & 100 & 19 (0) & 16.5 & 100 $\pm$ 0 & 12.4 $\pm$ 1.2 (0.0 $\pm$ 0.0) & 0.7 $\pm$ 0.6 & 0.0 $\pm$ 0.0 \\
(3, 3, 30)i & 96.7 & 34 (0) & 13.7 & 90.7 $\pm$ 3.3 & 24.3 $\pm$ 2.8 (16.0 $\pm$ 7.1) & 0.5 $\pm$ 0.3 & -1.3 $\pm$ 3.7 \\
(3, 3, 45)h & 77.8 & 49 (2) & 209.4 & 71.2 $\pm$ 3 & 29.2 $\pm$ 2.8 (12.5 $\pm$ 5.8) & 1.3 $\pm$ 0.8 & -5.2 $\pm$ 3.8 \\
\midrule

(6, 3, 30)l & 100 & 35 (0) & T/O & 100 $\pm$ 0 & 20.2 $\pm$ 1.2 (0.0 $\pm$ 0.0) & 3.3 $\pm$ 2.9  & 0.0 $\pm$ 0.0 \\
(6, 3, 60)i & 96.7 & 65 (0) & T/O & 94.4 $\pm$ 1.8 & 48.1 $\pm$ 4.3 (7.2 $\pm$ 5) & 7.8 $\pm$ 4.2 & -2.3 $\pm$ 3.6 \\
(6, 3, 90)h & 90 & 118 (3.4) & T/O & 80.4 $\pm$ 2.2 & 58.9 $\pm$ 3.9 (13.5 $\pm$ 6.9) & 20.8 $\pm$ 8.2 & -13.2 $\pm$ 3.5 \\

(6, 5, 30)l & 93.3 & 23 (0) & T/O & 93.3 $\pm$ 0 & 17.6 $\pm$ 1.4 (0.0 $\pm$ 0.0) & 5.8 $\pm$ 3.6 & 0.0 $\pm$ 0.0\\
(6, 5, 60)i & 98.3 & 62 (0) & T/O & 94.3 $\pm$ 2.4 & 48.4 $\pm$ 5.0 (6.7 $\pm$ 4.3) & 8.7 $\pm$ 4 & 1.1 $\pm$ 3.3 \\
(6, 5, 90)h & 75.6 & 81 (1.2) & T/O & 68.5 $\pm$ 1.6 & 54.0 $\pm$ 4.1 (17.6 $\pm$ 5.8) & 18.9 $\pm$ 10.8 & 1.7 $\pm$ 3.5\\
\midrule

(12, 3, 60)l & 98.3 & 47 (0) & T/O & 98.3 $\pm$ 0.0 & 36.3 $\pm$ 1.9 (0.0 $\pm$ 0.0) & 97.5 $\pm$ 76.3 & 0 $\pm$ 0 \\
(12, 3, 120)i & 34.2 & 42 (0) & T/O & 96.5 $\pm$ 1.2 & 93.8 $\pm$ 6.3 (11.0 $\pm$ 4.3) & 170.6 $\pm$ 74 & -17.2 $\pm$ 1.6\\
(12, 3, 180)h & - & - & M/O & 82.6 $\pm$ 1.7 & 123.6 $\pm$ 6.9 (14.2 $\pm$ 3.6)  & 443.9 $\pm$ 224.8 & - \\

(12, 5, 60)l & 98.3 & 56 (0) & T/O & 98.3 $\pm$ 0 & 32.8 $\pm$ 2.2 (0.0 $\pm$ 0.0) & 65.7 $\pm$ 55 & 0.2 $\pm$ 0.3\\
(12, 5, 120)i & 13.3 & 7 (0) & T/O & 95.5 $\pm$ 1.8 & 90.4 $\pm$ 8.8 (8.9 $\pm$ 3.6) & 156.2 $\pm$ 87.3 & -4.8 $\pm$ 3.8\\
(12, 5, 180)h & - & - & M/O & 85.1 $\pm$ 1.5 & 120.0 $\pm$ 9.8 (10.3 $\pm$ 3.4) & 455.3 $\pm$ 233.9 & - \\

(12, 7, 60)l & 98.3 & 41 (0) & T/O & 98.3 $\pm$ 0.0 & 40.3 $\pm$ 2.1 (0 $\pm$ 0) & 59.2 $\pm$ 42 & 2.0 $\pm$ 1.4 \\
(12, 7, 120)i & 13.3 & 14 (0) & T/O & 90.9 $\pm$ 1.3 & 93.9 $\pm$ 4.4 (14.9 $\pm$ 5.4) & 163.1 $\pm$ 90.3 & -18.8 $\pm$ 1.6 \\
(12, 7, 180)h & - & - & M/O & 80.6 $\pm$ 2.5 & 122.1 $\pm$ 8.6 (15.5 $\pm$ 3.8) & 491.8 $\pm$ 161.6 & - \\
\bottomrule

\end{tabular}}
\caption{Comparison between Gurobi 9.1.0 and our UAM-VNS on  (Aircraft, Vertiports, Requests) instances. \% Service is the percentage of served requests. \#Charge (\% Fast) is the number of charges and the percentage of fast ones. CPS Diff is the relative difference in Cost per Service, negative values indicating a lower CPS for the UAM-VNS. T/O (resp. M/O) indicates that the time (resp. memory) limit has been reached. }
    \label{tab:benchmark}
\end{table}

\subsubsection{Performance analysis}

Gurobi is able to prove optimality, within the time limit, only for the 3 smallest instances. For the 3 largest instances, routing 180 requests, it reaches memory out before producing any result. For all other intermediate instances, Gurobi is not able to prove optimality. We notice that for low-demand scenarios our UAM-VNS consistently finds solutions equivalent (in \% Service) to the Gurobi solutions. For the intermediate-demand scenarios, UAM-VNS achieved similar levels of service to Gurobi for small and average sized instances, i.e., aircraft fleets of sizes 3 and 6. 

As expected UAM-VNS is significantly faster than Gurobi. We can remark however that it relies on fast charging more often than Gubori while generally providing better CPS (negative difference), and similar service levels. Indeed, in several intermediate and high demand scenarios, with less than 6 aircraft, the UAM-VNS achieves a significantly lower CPS than Gurobi, this is explained in part by the lower number of served requests, allowing the UAM-VNS to select a subset of requests for which aircraft can be routed more cost-efficiently. Our method significantly outperforms Gurobi on larger instances (120 requests and 12 aircraft) both in term of service level and CPS ratio. 

\subsubsection{Discussion}

We proposed a VNS to solve the routing and recharging problem which can scale to large instances and yields solutions of good quality within acceptable time. Furthermore, the results presented in the previous paragraph lead to interesting perspectives.

First, we remarked that solutions obtained by UAM-VNS relied more on fast charges than those obtained by Gurobi, at the same time the UAM-VNS uses less charges in general. This can be explained by the design of neighborhood $N_3$ which performs charging moves. A local move consists in selecting the best charging mode and greedily optimizing the duration for a single charging operation. Therefore UAM-VNS is more likely to produce solutions with fewer and longer charging operations, rather than solutions with several shorter charging operations. This observation opens to future improvements of the charging neighborhood, for example by touching multiple charging operations in a single move. 

Second, we observed that in some cases the UAM-VNS achieved lower CPS than Gurobi while serving less requests. This suggests the possibility to select a subset of requests to serve based on CPS, for example level of service (number of requests served) could be maximized subject to a hard constraint on CPS.

\subsection{Online Pooling, Scheduling, Routing, and Recharging}\label{sec:online}

In this section we evaluate the capabilities of our algorithms in the following online setup. Assuming a set $D$ of previously booked demands and corresponding flight requests $R$, a user submits a new demand to the broker. Before presenting a bookable option to the user, the broker must determine if the demand can be served or not by updating the pooling and by probing the operator for routing feasibility. We conduct experiments to determine the response time obtained using UAM-Beam Search and UAM-VNS to solve this incremental problem.

\subsubsection{Specific settings}

We adopt one of the infrastructures used in the previous section, with 3 vertiports and 12 aircraft. To evaluate our \emph{accept-or-reject} capabilities, we consider a scenario where a set of $D$ demands is generated according to the process described in Section \ref{sec:settings-pooling}, these demands are spread uniformly randomly over all possible routes in the infrastructure. These were processed by our UAM-Beam Search, yielding a set of $R$ requests all served by a routing and recharging solution $S$, computed using UAM-VNS. Starting from this solution, \textit{one additional} demand is generated using the generation process of Section \ref{sec:settings-pooling} and its corresponding route is chosen randomly and uniformly among the possible routes. This additional demand is processed in two steps:

\begin{itemize}
    \item First, it is incorporated into the existing set of requests using UAM-Beam Search. Because this algorithm is by design incremental on the number of demands (see algorithm \ref{alg:beam}), adding one demand consists in adding one level to the existing search tree. When this is done, an updated set of requests containing $R|D+1$ requests is available. At this stage, it is important to note that the new set of requests contains at most one additional request and that most requests are not affected by the addition of the new demand.
    \item Second, the updated set of requests is fed to our UAM-VNS to find a new solution. Because the new set of requests is by construction close to the initial set, the previous solution $S$ can be used as starting point for the UAM-VNS. Since the goal is to either accept of reject the new demand, the anytime property of the algorithm can be leveraged and we can stop the search as soon as a solution serving all requests is found. If no solution serving all requests is found, then the additional demand is not accepted.
\end{itemize}

In this experiment, we are interested in measuring acceptance status (i.e. indicating whether the additional demand is accepted or not) after different computation times, i.e. after 5, 10, 20, 30 seconds, and when the UAM-VNS stops by itself due to timeout or stagnation. For initial set of demands of several sizes, we generate an additional demand and execute the two steps described above 20 times. We report the percentage of acceptance after different running times in Table \ref{tab:online}.

\begin{table}[htbp]
\centering
\resizebox{\textwidth}{!}{
 \begin{tabular}{lc|cccccc}
\toprule

D & R & $R | D+1$ & \%Accept 5'' & \%Accept 10'' & \%Accept 20''  & \%Accept 30'' & \%Feasibility  \\

\hline
100  & 48 & 48.4 & 100 & 100 & 100 & 100 & 100   \\
120  & 54 & 54.3 & 100 & 100 & 100 & 100 & 100  \\
140  & 61 & 61.3 & 90 & 100 & 100 & 100 & 100   \\
160  & 62 & 62.3 & 100 & 100 & 100 & 100 & 100  \\
180  & 72 & 72.5 & 65 & 100 & 100 & 100 & 100  \\
200  & 78 & 78.5 & 35 & 65 & 75 & 95 & 100 \\
220  & 82 & 82.4 & 5 & 25 & 50 & 75 & 95 \\
\end{tabular}}
\caption{Demand acceptance in function of computation time (UAM Beam Search \& UAM-VNS).}
    \label{tab:online}
\end{table}

\subsubsection{Discussion}

One important feature of the acceptance process described in this section is that we never accept a demand that cannot be accommodated. We might however reject a demand when it would be possible to accommodate it using more computation time. The results indicate that allowing a time budget of 30 seconds to decide whether to accept a new demand or not resulted in wrongfully rejecting a demands 5\% to 25\% of the time only in the most constrained cases with 200 and 220 initial demands. When the number of demands is between 100 and 180, an additional demand can be accepted within 10 seconds 100\% of the time, and in under 5 seconds between 65\% and 100\% of the time.
\newline 

These results show that it is possible to perform the demand acceptance check sufficiently apace for online human interaction. 
One key aspect is that accepted demands act as engagements for the broker/operator(s). They will be performed as agreed, even if their cost might change due to further engagements. We will see in  the Conclusion section that demands can also be accepted depending on their potential future benefits, e.g., depending on how each new demand might be eventually pooled with future demands. This offers a perspective to mix prediction and optimization for both the broker and the operators.

\section{Conclusion and Perspectives}

We have presented the UAM on-demand air-taxi problem. A new transportation service which pools passengers into short urban flights connecting city-centers, important suburbs and airports.
We envisioned two categories of actors for this new service. Brokers market flights to passengers, and flight operators perform these flights through their eVTOL fleet.
Acknowledging the underlying conflicting incentives for optimization, maximizing pooling against operating as many flights as possible, we have decomposed this problem into two parts. First, the pooling and scheduling of passengers into joint flights. Second, the routing and recharging of these flights through an UAM infrastructure.

Contrary to traditional airlines, restricted eVTOLs size and design will not allow different space or catering experiences for passengers.
Therefore, we proposed a passenger-centric formulation for the broker's Pooling and Scheduling problem where the quality of service is maximized by minimizing passenger waiting times for on-boarding. Our approach allows to differentiate service levels through multiple passenger classes, e.g., Regular, and Premium.
The second part, Routing and Recharging, was formulated such that routing and energy recharging decisions are taken jointly, including the possibility of using several charging modes to efficiently execute flights while favoring long term battery health. 
For each sub-problem, we proposed dedicated online optimization algorithms which provide good solutions and scale better than the latest version of a state-of-the-art commercial solver, Gurobi.

Our experiments and detailed analysis provided meaningful insights into the design space of the Pooling and Scheduling problem. We highlighted how quality of service differs between passenger classes as problem parameters vary. Our analysis also revealed potential fairness issues and we provided both explanations and solutions to this important problem. 
In the case of more than 2 classes, we expect that our analysis and the fairness issues  would still hold and that these could be exacerbated when comparing classes that are far from each other on the classes spectrum. 
We believe that these insights could help designing future air-taxi services and determine appropriate service quality and pricing strategies. 
Furthermore, the time-constrained experiments presented in section \ref{sec:online} indicate that, thanks to the good performance of our algorithms, the air-taxi option could be embedded smoothly into existing transportation platforms where end-users would be presented with the air-taxi option only when available.  \\


There are several perspectives for this work. First, at the technical and algorithmic level, we will work on improving the performance and the quality of the solutions returned by the UAM-VNS algorithm. As identified in experiments, designing a neighborhood for electric recharging actions that modifies several operations simultaneously could yield better solutions, using less fast charges.
Another very promising line of work is to use predictive models.
For instance, when there are different ways to insert deadheads while building a routing solution, an operator-request prediction model could be used to place a deadhead where new flight requests are most likely to occur later.
Similarly, a broker-demand prediction model could be used to preemptively reserve flights to the operator.
To optimize against uncertainty, we could leverage robust and stochastic optimization techniques able to tolerate perturbations and to adapt solutions through recourse operators for example.\\

Beyond algorithmic extensions, our passenger-centric approach could be augmented to include a deeper level of personalization. For example, in order to guarantee a good user-experience, and considering the fairness issues we identified, it could be interesting to use a record of QoS levels experienced by each passenger in order to mitigate fairness issues over several trips. For instance, this could be done by not always pooling the same regular passenger with premium passengers. 
These considerations centered around user experience and QoS will help differentiate UAM air-taxi brokers from other transportation services. \\

This paper introduces a rich context for urban air-taxi services and provides the first algorithmic methods to operate these at scale. Our results and analysis open to very interesting perspectives and we hope this will inspire future research in this domain. Furthermore, our study can help actual air-taxi developers \cite{Volocopter2019,ehang2020,whitepapercommunityelevate2020,UberWP2016} to design products that maximize quality of service while treating passengers fairly.

\bibliographystyle{plain}
\bibliography{references}
\end{document}